\documentclass[12pt]{amsart}
\usepackage{graphicx,amssymb}
\newtheorem{thm}{Theorem}[section]
\newtheorem{cor}[thm]{Corollary}
\newtheorem{lem}[thm]{Lemma}
\newtheorem{prop}[thm]{Proposition}
\theoremstyle{definition}
\newtheorem{defn}[thm]{Definition}
\theoremstyle{remark}
\newtheorem{rem}[thm]{Remark}
\numberwithin{equation}{section}

\begin{document}

\title[Smooth $K$-theory ]{Smooth $K$-theory of locally convex algebras}%
\author{H. Inassaridze and T. Kandelaki}%
\address{M.Alexidze Str.1, Tbilisi 0193, Georgia}%
\email{hvedri@rmi.acnet.ge  and  kandel@rmi.acnet.ge}%

\thanks{}%
\subjclass{19D50, 46L80, 19D55, 19L99, 18G55}%
\keywords{ smooth map, approximate unit, locally convex algebra,
smooth homotopy,
Fr\'echet algebra, algebraic $K$-functor.}%

\begin{abstract}Smooth $K$-functors are introduced and the smooth
$K$-theory of locally convex algebras is developed. It is proved
that the algebraic and smooth $K$-functors are isomorphic on the
category of quasi $\hat\otimes$-stable real (or complex) Fr\'echet
algebras.
\end{abstract}
\maketitle
\section*{Introduction}

We develop a new topological $K$-theory for arbitrary locally
convex $k$-algebras called smooth $K$-theory and constructed by
using smooth maps.

The category of locally convex $k$-algebras is a wide class of
topological $k$-algebras, containing Fr\'echet $k$-algebras,
Michael's $k$-algebras which are isomorphic to projective limits
of Banach $k$-algebras \cite{Mich}, many important examples of
differential operator algebras and differential forms, closely
related to noncommutative geometry.

The definition of smooth $K$-functors is motivated by our purpose
to extend Karoubi's Conjecture on the isomorphism of algebraic and
topological $K$-functors \cite{Karo} to a wide class of
topological $k$-algebras containing those for which Karoubi's
Conjecture was already confirmed, namely the category of
$C^*$-algebras \cite{SuWo}, the category of generalized operator
algebras and their polynomial extensions \cite{InKa}, by using the
important notion of stability of these algebras. For more details
see the summarizing article \cite{Rose} about the relationship
between algebraic and topological K-theory for Banach algebras and
C*-algebras. The real case will be also treated. In connection
with this problem the smooth $K$-functors appear in a natural way
and play a fundamental role in establishing our main result
(Theorem 4.4) confirming what we call the Smooth Karoubi's
Conjecture:

{\em The algebraic and smooth $K$-functors are isomorphic on the
category of quasi $\hat{\otimes}$-stable real (or complex)
Fr\'echet algebras.}

Consequently the use of smooth maps leads to an unexpected
relationship between algebraic $K$-theory and topologically
defined smooth $K$-theory via Grothendieck projective tensor
product of locally convex algebras.

The proof of this theorem is spread over Sections 1,2,3 and 4.

It should be noted that Phillips \cite{Phil} has defined a
K-theory for complex Frechet algebras that are locally
multiplicatively convex, extending the topological K-theory of
complex Banach algebras and showing that many important
K-theoretical properties are preserved such as homotopy
invariance, exactness, periodicity and stability properties. A
bivariant $K$-functor on the category of complex locally
multiplicatively convex algebras has been constructed by Cuntz
\cite{Cunt1} defining a bivariant multiplicative character to
bivariant periodic cyclic cohomology. When the first variable is
trivial, Cuntz's construction provides a K-theory isomorphic to
Phillips' K-theory. Both K-theories of complex locally
multiplicatively convex Fr\'echet algebras are constructed in
different context.

To prove the Smooth Karoubi's Conjecture we have restricted
ourselves by considering Fr\'echet $k$-algebras which topologies
are determined by countably many seminorms. We don't assume that
these seminorms possess the multiplicative property. The proof is
essentially based on the countability property of determining
seminorms. It seems that our main theorem is not valid in the
category of arbitrary (non-unital) locally convex $k$-algebras.
This is due to the fact that algebraic $K$-functors in general
case are not smooth homotopy functors for quasi
$\hat\otimes$-stable locally convex $k$-algebras having a bounded
approximate unit, since they don't possess probably the
Suslin-Wodzicki's TF-property of rings and the excision property.

In the case of Fr\'echet $k$-algebras, in particular for Banach
$k$-algebras, the $k$-algebra $A^I$ of all continuous maps
$I=[0,1]\rightarrow A$ into a Fr\'echet $k$-algebra $A$ is not
compatible with the Grothendieck projective tensor product and we
don't have the following needed isomorphisms:
\begin{equation}\label{inteqa}
(A\hat{\otimes}k)^I\approx A\hat{\otimes} k^I\approx
A^I\hat{\otimes}k.
\end{equation}
That is the reason why the usual construction of topological
$K$-functors seems to be difficult to be used to confirm Karoubi's
Conjecture if the stability is expressed in terms of the
Grothendieck projective tensor product. On the other hand, the
$k$-algebra $A^{\infty(I)}$ of all smooth maps $I\rightarrow A$
possesses the property (\ref{inteqa}) and many other functorial
properties helping us to prove our main theorem.

In Section 1 smooth $K$-functors are introduced and investigated
in the category of locally convex $k$-algebras. Their main
homological properties are established such as exactness, smooth
homotopy property, compatibility with finite products,
relationship with algebraic and topological $K$-theories in terms
of sufficient conditions implying isomorphisms of smooth
$K$-theory with these $K$-theories. It is shown that the smooth
$K$-theory agrees to Phillips' $K$-theory.

Section 2 deals with the Cohen-Hewitt factorization theorem. It is
well-known that the excision property in algebraic $K$-theory is
closely related with the triple factorization ($TF$-) property of
rings \cite{SuWo}. The $TF$-property for Banach algebras is
established in \cite{SuWo} by using a theorem of Hewitt
\cite{Hewi} motivated by Cohen's result \cite{Cohe}. According to
Cohen-Hewitt theorem any element of a Banach module can be
factorized into a product of two elements, when the Banach algebra
possesses a bounded approximate unit which is an approximate unit
of the module. In \cite{Wodz} Wodzicki showed that every
multiplicatively convex Fr\'echet $k$-algebra with uniformly
bounded approximate unit is $H$-unital as a consequence of the
Cohen-Hewitt factorization theorem extended to the case of
multiplicatively convex Fr\'echet algebras. In this section the
Cohen-Hewitt factorization theorem will be proved for Fr\'echet
modules over a Banach $k$-algebra which have a bounded approximate
unit.

In Section 3 we give an account on Higson's homotopy invariance
theorem. This important theorem plays an essential role for
setting the homotopy invariance of functors (for example, using
this theorem Karoubi's Conjecture on the isomorphism of algebraic
and topological $K$-theories was confirmed on the category of
stable $C^{*}$-algebras \cite{SuWo}). The natural question arises
whether this theorem is true or not for real $C^{*}$-algebras,
since we could not extend the proof of Lemma 3.1.2 and Theorem
3.1.1 of \cite{Higs} to the case of real $C^{*}$-algebras. In this
section we confirm that Higson's homotopy invariance theorem holds
also for functors defined on the category of real
$C^{*}$-algebras.

Section 4 is devoted to the Smooth Karoubi's Conjecture. We prove
that Fr\'echet $k$-algebras with bounded approximate unit have the
triple factorization property and therefore they possess the
excision property in algebraic K-theory. Using Higson's homotopy
invariance theorem we conclude that the functors
$K_{n}(-{\hat\otimes {\mathcal{K}}})$, $n\geq 1$, are smooth
homotopy functors implying the confirmation of the Smooth
Karoubi's Conjecture for real and complex cases.

 We close this introduction with a list of
terminological and notational conventions used in the present
article:

1) $k$ denotes the field of real or complex numbers;

2) algebras are always associative and are not assumed in general
to possess unit;

3) $A^{+}= A+k$ denotes the $k$-algebra obtained by adjoining unit
to a $k$-algebra $A$;

4) $\mathcal{K}$ denotes the $C^{*}$-algebra of compact operators
on the standard infinite-dimensional separable Hilbert space;

5) $\hat{\otimes}$ denotes the Grothendieck projective tensor
product of locally convex $k$-algebras;

6) $\mathcal{G}r$ and $\mathcal{A}b$ denote respectively the
category of groups and the category of abelian groups;

7) $Ob\mathcal{A}$ denotes the class of objects of the category
$\mathcal{A}$;

8) the determining seminorms of Fr\'echet $k$-algebras are not
assumed to possess the multiplicative property;

9) without loss of generality a locally convex $k$-algebra with
bounded approximate unit always means with left bounded
approximate unit and having the TF-property always means having
the left triple factorization property.

\

\section{Smooth $K$-theory}

The smooth $K$-functors $K^{sm}_n,\;\;n\geq 0$, will be defined on
the category $\mathcal{A}$ of locally convex $k$-algebras and
their continuous $k$-homomorphisms.

By a locally convex $k$-algebra we mean a $k$-algebra $A$ equipped
with a complete locally convex topology such that the
multiplication map $A\times A\rightarrow A$ is jointly continuous.
The category $\mathcal{A}$ is close under the Grothendieck
projective tensor product \cite{Trev}. We recall its definition.

Let $A$ and $B$ be two locally convex $k$-algebras. Their
projective tensor product $A\hat{\otimes}B$ is given by a family
$I$ of seminorms $\mu\otimes\nu$ on $A\otimes B$,
\begin{equation}
\label{smeqa} (\mu\otimes \nu)(x)=\inf \bigg( \sum
_{k=1}^n\mu(a_k)\nu (b_k)\bigg),
\end{equation}
where infimum is taken over all representations of $x$ of the form
$$
x=\sum _{k=1}^n a_k\otimes b_k,\;\;\;\;\;a_k\in A,\;\;b_k\in B,
$$
the locally convex $k$-algebra $A\hat{\otimes}B$ being the
completion of $A\otimes B$ with respect to the family
$I=\{\mu{\otimes}\nu\}$ of seminorms.

The projective tensor product possesses exactness property with
respect to proper short exact sequences in the category of locally
convex $k$-algebras that we shall need in Section 4.

 Let
$$
0\rightarrow I\xrightarrow {f}B\xrightarrow {g}A\rightarrow 0
$$
be a sequence of morphisms in the category $\mathcal{LC}$ of
locally convex linear topological spaces and continuous linear
maps. It will be said that this sequence is a {\em proper} exact
sequence if $f$ is a homeomorphism of $I$ on $Imf$ and $g$ is an
open surjective map. It will be said {\em proper split} exact if
$g$ has a right inverse in $\mathcal{LC}$. A short exact sequence
in the category $\mathcal{A}$ is said {\em proper} exact sequence
if it is proper split exact in the category $\mathcal{LC}$.
\begin{lem} Let
$$
0\rightarrow I\xrightarrow {f}B\xrightarrow {g}A\rightarrow 0
$$
be a proper exact sequence of locally convex $k$-algebras and $H$
be a locally convex $k$-algebra. Then the sequence
$$
0\rightarrow
H\hat{\otimes}I\xrightarrow{H\hat{\otimes}f}H\hat{\otimes}B\xrightarrow
{H\hat{\otimes}g}H\hat{\otimes}A\rightarrow 0
$$
is a proper exact sequence.
\end{lem}
\begin{proof} Let $l:A\rightarrow B$ be a continuous linear map
such that $gl= 1_{A}$. It is clear that $lg:B\rightarrow B$ is a
continuous projection and $Im(lg)$ is a closed linear subspace of
$B$ which is isomorphic to $A$. It is easily checked that the
natural continuous linear map $I\times A\rightarrow B$ given by
$(x,y)\mapsto f(x)+l(y)$ is an isomorphism of linear topological
spaces. The converse map is defined by $z\mapsto
(f^{-1}(z-lg(z)),g(z))$. One has the following commutative diagram
\begin{equation}
\label{appeqc}
\begin{array}{ccccccccc}
 0 & \rightarrow & I \hat{\otimes}H & \xrightarrow{(1_{I \hat{\otimes}H},0)} &
  (I \hat{\otimes}H)\times (A \hat{\otimes}H)
  & \xrightarrow{pr_2} & A \hat{\otimes}H & \rightarrow & 0 \\
  &  & \| &  & \downarrow &  & \| &  &  \\
0 & \rightarrow & I \hat{\otimes}H & \xrightarrow{f\hat{\otimes}
1_H} & B \hat{\otimes}H
  & \xrightarrow{g\hat{\otimes}1_H} &
  A \hat{\otimes}H & \rightarrow & 0 , \\
\end{array}
\end{equation}
where the top row is the trivial proper split exact sequence of
locally convex $k$-algebras and the middle vertical arrow is the
composite of the natural isomorphisms $(I\hat{\otimes}H)\times
(A\hat{\otimes}H)\xrightarrow {\approx} (I\times
A)\hat{\otimes}H\xrightarrow {\approx}{B\hat{\otimes}H.}$
Therefore the bottom row is a proper split short exact sequence of
locally convex $k$-algebras.
\end{proof}

It is said that a locally convex $k$-algebra $A$ possesses  a left
bounded approximate unit if there exists a bounded direct set
$\{e_{\lambda}\}$ of elements in $A$ such that
$$
\underrightarrow{\lim}_i\parallel e_{\lambda}a-a\parallel_{i}=0
$$
for all $a\in A$ and $i\in I$. The right bounded approximate unit
in $A$ is defined similarly. It is easily checked that the
category of locally convex $k$-algebras with left bounded
approximate unit is closed under the projective tensor product.

 A Fr\'echet $k$-algebra is a locally convex complete $k$-algebra
such that its topology is given by a countable family of
seminorms.

The category $\mathcal{B}$ of Fr\'echet $k$-algebras is a full
subcategory of the category $\mathcal{A}$ and is closed under the
Grothendieck tensor product.

A continuous map $f:I\rightarrow A $ with values in a locally
convex $k$-algebra $A$ is called smooth, if it has all derivatives
$f^{(1)},f^{(2)},...,f^{(n)},....$ The space $A^{\infty (I)}$ of
all smooth maps from $I$ to $A$ is also a locally convex
$k$-algebra. The seminorms $\mu _{s,r}$ on $A^{\infty (I)}$ are
given by
\begin{equation}
\mu _{s,r}(f)=\sup _t \{\sum _{k=1}^r ||f^{(k)}(t)||_s,\;\;\;,t\in
I\},
\end{equation}
where $\{|||-|||_s\}_{s\in S}$ is the set of the determining
seminorms on $A$. If $A$ is a Michael $k$-algebra or a Fr\'echet
$k$-algebra, then $A^{\infty (I)}$ is also a Michael $k$-algebra
or a Fr\'echet $k$-algebra respectively.

Any continuous $k$-homomorphism $\varphi :A\rightarrow A'$ of
locally convex $k$-algebras induces in a natural way a continuous
$k$-homomorphism $\varphi ^{\infty (I)}:A^{\infty (I)}\rightarrow
A'^{\infty (I)}$.

Now on the category $\mathcal{A}$ we will define the smooth path
cotriple $\mathcal{I}$.

The evolution maps at $t=0$ and $t=1$,
$$
\varepsilon_i :A^{\infty (I)}\rightarrow
A,\;\;i=0,1,\;\;\varepsilon_0(f)=f(0),\;\; \varepsilon _1(f)=f(1),
$$
will play an important role. Denote by $\mathcal{I}(A)$ the kernel
of $\varepsilon _0$ and by $\tau _A:\mathcal{I}(A)\rightarrow A$
the restriction of $\varepsilon _1$ on $\mathcal{I}(A)$. There is
a natural continuous $k$-homomorphism $\delta _A:
\mathcal{I}(A)\rightarrow \mathcal{I}^2(A)=\mathcal{I}
(\mathcal{I} (A))$ sending $f\in \mathcal{I}(A)$ to $\delta
_A(f)(s,t)=f(st)$. It is easily checked that $\delta _A(f)$ is a
smooth map.

By taking $\textbf{I}=(\mathcal{I},\tau , \delta)$ one gets a
cotriple on the category $\mathcal{A}$ which will be called the
smooth path cotriple. The cotriple $\mathbf{I}$ induces the
augmented simplicial locally convex $k$-algebra
\begin{equation}
\label{smeqb}
  \mathcal{I}_*^+(A)=\mathcal{I}_*(A)\rightarrow A,
\end{equation}
where $\mathcal{I}_0(A)=\mathcal{I}(A)$,\;
$\mathcal{I}_n(A)=\mathcal{I}(\mathcal{I}_{n-1}(A)),\;\;n\geq 1,$
\;$\delta _i^n=\mathcal{I}^i\tau \mathcal{I}^{n-i}$,\;
$s_i=\mathcal{I}^i\delta \mathcal{I}^{n-i}$, \;$n\geq 1$,\; $0\leq
i\leq n$.

Applying the general linear group functor one obtains an augmented
simplicial group
$$
GL(\mathcal{I}_*^+(A))=GL(\mathcal{I}_*(A))\rightarrow GL(A).
$$
\begin{defn}
\label{smdfa} We define the smooth $K$-functors
$K^{sm}_n,\;\;n\geq 0$, by setting
$$
K^{sm}_n(A)=
\left\{%
\begin{array}{ll}
   \pi _{n-2}GL(\mathcal{I}_*(A)) & \hbox{for $n\geq 3$,} \\
    K_0(A), & \hbox{$n=0$,} \\
\end{array}%
\right.
$$
and $K^{sm}_1(A),\;\;\;K^{sm}_2(A)$ make exact the following
sequence
$$
0\rightarrow K_2^{sm}(A)\rightarrow
\pi_0GL(\mathcal{I}_*(A))\xrightarrow{\gamma _0} GL(A)\rightarrow
K_1^{sm}(A)\rightarrow 0,
$$
where $\gamma _0$ is induced by $GL(\tau_A)$, or equivalently
$$
K_n^{sm}(A)=\pi _{n-2}GL(\mathcal{I}_*^+(A)),\;\;\;n\geq 1.
$$
\end{defn}
From the definition of smooth $K$-functors follows immediately
$$K_n^{sm}(\mathcal{I}(A))= 0$$ for all $n\geq 1$ and any locally
convex $k$-algebras, since the augmented simplicial group
$GL\mathcal{I}_*^+(\mathcal{I}(A))$ is right contractible and
therefore aspherical according to Lemma 1.1 and Lemma 1.2
\cite{Swa}.

It will be shown that all smooth $K$ -groups are abelian groups.
This is obvious for $n=0$ and $n\geq 3$. First we will proof that
$K^{sm}_1(A)$ is a group or equivalently that
$\mathrm{Im}\;GL(\tau _A)$ is a normal subgroup of $GL(A)$. The
exact sequence
$$
0\rightarrow \mathcal{I}(A)\xrightarrow{\sigma} A^{\infty
(I)}\xrightarrow{\varepsilon_0} A\rightarrow 0
$$
yields the short exact sequence
$$
0\rightarrow GL(\mathcal{I}(A))\xrightarrow{GL(\sigma )}GL(
A^{\infty (I)})\xrightarrow{GL(\varepsilon _0)} GL(A)\rightarrow
0,
$$
since $GL(\varepsilon _0)$ is splitting. Thus $GL(\mathcal{I}(A))$
is a normal subgroup of $GL( A^{\infty (I)})$. The assertion
follows from the fact that $GL(\varepsilon _1)$ is surjective.

Consider the free cotriple $\textbf{F}=(\mathcal{F},\tau ,\delta)$
in the category of rings. Then for any locally convex $k$-algebra
$A$ one has a homomorphism $\alpha _0:\mathcal{F}(A)\rightarrow
\mathcal{I}(A)$ given by $|a|\mapsto (t\mapsto at)$, $a\in
A,\;\;t\in I$. It is clear that the map $t\mapsto at,\;\;t\in I,$
is a smooth map. We have also the inclusion $\beta
_0:\mathcal{I}(A)\rightarrow \mathcal{J}(A)$, where
$\mathbf{J}=(\mathcal{J},\tau ,\delta)$ is the continuous path
cotriple in the category $\mathcal{A}$. As a consequence one has
morphisms of cotriples
$$
\mathbf{F}\xrightarrow{\alpha} \textbf{I} \xrightarrow{\beta}
\mathbf{J}
$$
such that the composition $\beta\alpha :\mathbf{F}\rightarrow
\mathbf{J}$ is the well-known morphism from the free cotriple to
the continuous path cotriple.

For any locally convex $k$-algebra $A$ the groups
$$
K^{top}_{n}(A)=\pi_{n-2}GL(L^{+}_{*}(A))\;,n\geq 1,
$$
and $K^{top}_{0}(A)=K_{0}(A)$ are called the topological
$K$-groups of $A$. In fact these topological $K$-groups have been
already defined by Swan for any real topological algebra [18].
When $A$ is a Banach $k$-algebra we recover the well-known
topological $K$-groups of $A$ \cite{Karo}.

 The morphisms $\alpha $ and $\beta $
induce respectively functorial homomorphisms
$$
K^{S}_n(A)\xrightarrow{\alpha ^*_n} K^{sm}_n(A)\xrightarrow{\beta
^*_n}K^{top}_n(A)
$$
for $n\geq 0$, where $K^{S}_*$ are Swan's $K$-functors which are
isomorphic to Quillen's $K$-functors $K_{*}$.

It is clear that the homomorphisms
$$
\alpha ^*_1:K^{S}_1(A)\rightarrow
K_1^{sm}(A),\;\;\;\;\;\;\;\;\beta ^*_1:K_1^{sm}(A)\rightarrow
K_1^{top}(A)
$$
are surjective and therefore $K_1^{sm}(A)$ is an abelian group.
\begin{defn}
\label{smdfb} A continuous $k$-homomorphism $f:A\rightarrow A''$
of locally convex $k$-algebras is called $GL$-fibration with
respect to the smooth path cotriple $\mathbf{I}$, if the induced
homomorphisms
$$
GL(\mathcal{I}_i(A))\rightarrow GL(\mathcal{I}_i(A''))
$$
are surjective for all $i\geq 0$.
\end{defn}
It is obvious that any splitting continuous $k$-homomorphism
$f:A\rightarrow A''$  is a $GL$-fibration with respect to
$\mathbf{I}$.

 Thus the sequence
$$
0\rightarrow A\xrightarrow{i} A^{+}\xrightarrow{p} k\rightarrow 0
$$
is a $GL$-fibration with respect to the smooth path cotriple
$\mathcal{I}$, where $i(a)=(a,0)$ and $p(a,t)=t,\; a\in A,\; t\in
k$ and this short exact sequence yields the splitting exact
sequence
$$
0\rightarrow K^{sm}_{n}(A)\rightarrow K^{sm}_{n}(A^{+})\rightarrow
K^{sm}_{n}(k)\rightarrow 0
$$
for $n\geq 0$.

\begin{prop}
\label{smprpa} Any short exact sequence of locally convex
$k$-algebras
$$
0\rightarrow A'\xrightarrow{\sigma} A\xrightarrow{\eta}
A''\rightarrow 0\;,
$$
where $\eta$ is a $GL$-fibration with respect to the smooth path
cotriple $\mathbf{I}$, induces a long exact sequence of smooth
$K$-functors
\begin{multline}
\label{smmlta} ...\rightarrow K^{sm}_{n+1}(A'')\rightarrow
K^{sm}_{n}(A')\rightarrow K^{sm}_{n}(A)\rightarrow
K^{sm}_{n}(A'')\rightarrow K^{sm}_{n-1}(A')\rightarrow ...\\
...\rightarrow K^{sm}_{2}(A'')\rightarrow
K^{sm}_{1}(A')\rightarrow K^{sm}_{1}(A)\rightarrow
K^{sm}_{1}(A'')\rightarrow \\\rightarrow K_{0}(A')\rightarrow
K_{0}(A)\rightarrow K_{0}(A'').
\end{multline}
\end{prop}

\begin{proof}
 The given short exact sequence of locally convex $k$-algebras
 yields the short exact sequence of augmented simplicial groups
$$
0\rightarrow GL(\mathcal{I}^+_*(A'))\rightarrow
GL(\mathcal{I}^+_*(A))\rightarrow
GL(\mathcal{I}^+_*(A''))\rightarrow 0
$$
which implies the required long exact sequence of smooth
$K$-functors ending with $K_1^{sm}(A)$. It remains to show the
exactness of
$$
K^{sm}_{1}(A)\rightarrow K^{sm}_{1}(A'')\xrightarrow{\partial _1}
K_{0}(A')\rightarrow K_{0}(A)\rightarrow K_{0}(A'').
$$
The following commutative diagram
$$
\begin{array}{ccccccccccc}
    && GL(\mathcal{I}(A)) & \rightarrow & GL(\mathcal{I}(A'')) & \rightarrow & 0 &&& &  \\
     &  & \downarrow &  &  \downarrow &  &  & & & &  \\
    ... & \rightarrow & K_1(A) & \rightarrow & K_1(A'') & \rightarrow
   & K_0(A')& \rightarrow & K_0(A) &  &  \\
   &  & \downarrow &  & \downarrow &  &  &  &  &  &  \\
 K_1^{sm}(A') & \rightarrow &  K_1^{sm}(A) & \rightarrow &
K_1^{sm}(A'') &
  &  &  &  &  &  \\
    &  & \downarrow &  & \downarrow &  & &  &  &  &  \\
   &  & 0 &  & 0 &  &  &  &  &  &  \\
\end{array}
$$
with exact rows and columns defines $\partial _1$ in a natural way
and implies the exactness of the remaining sequence.
\end{proof}

An assertion similar to Proposition 1.4 holds also for the
topological $K$-functors $K_{n}^{top}$, $n\geq o$, with respect to
the continuous path cotriple $J$.

It is easy to show that the sequence
\begin{equation}
\label{smeqc} 0\rightarrow \Omega_{sm} (A)\rightarrow
\mathcal{I}(A)\xrightarrow{\tau _A} A\rightarrow 0
\end{equation}
is a $GL$-fibration with respect to $\mathbf{I}$, where
$\Omega_{sm} (A)=\mathrm{Ker}\;\tau _A$, the homomorphisms
$GL\mathcal{I}_i(A)\rightarrow GL\mathcal{I}_{i+1}(A)$ given by
$\mathcal{I}_i(\tau _A)$, $\mathcal{I}_0(\delta _A)=\delta _A$,
$\mathcal{I}_i(\delta _A)=\mathcal{I}(\mathcal{I}_{i-1}(\delta
_A))$, for $i>0$, being the splitting homomorphisms for
$GL\mathcal{I}_i(\tau _A)$, $i\geq 0$.

Applying Proposition \ref{smprpa} to the sequence \ref{smeqc} one
obtains the isomorphism $K_2^{sm}(A)\approx K_1^{sm}(\Omega_{sm}
(A))$ showing that $K_2^{sm}(A)$ is an abelian group.

\begin{defn}
\label{smdfc} Two continuous $k$-homomorphisms $f,\;g:A\rightarrow
B$ are called smooth homotopic if there exists a continuous
$k$-homomorphism $h:A\rightarrow B^{\infty (I)}$ such that
$\varepsilon _0h=\varepsilon _1h$, which is called the smooth
homotopy between $f$ and $g$.
\end{defn}
\begin{defn}
\label{smdfd} A functor $T:\mathcal{A}\rightarrow \mathcal{G}r$ is
called a smooth homotopy functor if $T(f)=T(g)$ for smooth
homotopic continuous $k$-homomorphisms $f$ and $g$.
\end{defn}
It is obvious that $T$ is a smooth homotopy functor if
$T(\varepsilon _0)=T(\varepsilon _1)$.

The topological $K$-functors $K^{top}_{n}$ satisfy the condition
$K^{top}_{n}(\varepsilon_{0})=K^{top}_{n}(\varepsilon_{1})$,
$n\geq 1$, and therefore are smooth $K$-functors for all $n\geq
1$. It is well-known that on the subcategory of Banach
$k$-algebras the Grothendieck $K$-functor $K_{0}$ is a homotopy
functor implying by the same reason that $K_{0}$ is a smooth
homotopy $K$-functor on the category of Banach $k$-algebras.
\begin{prop}
\label{smprpb} A functor $T:\mathcal{A}\rightarrow \mathcal{G}r$
is a smooth functor if and only if the inclusion $i:A\rightarrow
A^{\infty (I)}$ induces an isomorphism
$T(i):T(A)\xrightarrow{\approx}T(A^{\infty (I)})$ for all $A\in
Ob\mathcal{A}$.
\end{prop}

\begin{proof}
If $T(i)$ is an isomorphism, then the equality $T(\varepsilon
_0)T(i)=T(\varepsilon _1)T(i)$ implies $T(\varepsilon
_0)=T(\varepsilon _1)$ and therefore $T$ is a smooth homotopy
functor. Conversely, let $T$ be a smooth homotopy functor. Then
the continuous $k$-homomorphism $i\varepsilon _0:A^{\infty
(I)}\rightarrow A^{\infty (I)} $ is smooth homotopic to the
identity map, where the map $A^{\infty (I)}\rightarrow
\big(A^{\infty (I)}\big)^{\infty (I)}$, sending the smooth map
$f:I\rightarrow A$ to the smooth map $\varphi :I\rightarrow
A^{\infty (I)}$, \;$\varphi (t)(x)=f(tx)$, $i,x\in I$, provides a
smooth homotopy between $i\varepsilon _0$ and $1_{A^{\infty
(I)}}$. Thus $T(i)T(\varepsilon _0)=id=T(\varepsilon _0)T(i)$. The
proof is complete.
\end{proof}
\begin{prop}
\label{smprpc} The smooth $K$-functors $K^{sm}_n,\;n\geq 1$, are
smooth homotopy functors.
\end{prop}
\begin{proof}
The exact sequence \ref{smmlta} applied to the smooth
$GL$-fibration
$$
0\rightarrow \mathcal{I}(A)\rightarrow A^{\infty (I)}
\xrightarrow{\varepsilon _0}A\rightarrow 0
$$
yields the short exact sequence
$$
0\rightarrow K^{sm}_n(\mathcal{I}(A))\rightarrow
K^{sm}_n(A^{\infty (I)}) \xrightarrow{\varepsilon
_0}K^{sm}_n(A)\rightarrow 0
$$
for $n\geq 1$. It remains to recall that
$K^{sm}_n(\mathcal{I}(A))=0$, $n\geq 1$, and one obtains the
isomorphism $K^{sm}_n((A))\xrightarrow{\approx}K^{sm}_n(A^{\infty
(I)})$, $n\geq 1$, showing by Proposition 1.7 that the
$K$-functors $K^{sm}_n$, $n\geq 1$, are smooth homotopy functors.
\end{proof}

The smooth homotopization of a functor $T:\mathcal{A}\rightarrow
\mathcal{G}r$ is a functor $h^{sm}T:\mathcal{A}\rightarrow
\mathcal{G}r$ given by
\begin{equation}
h^{sm}T(A)=\mathrm{Coker}(T(A^{\infty (I)})\rightrightarrows T(A))
\end{equation}
for any $A\in \mathrm{Ob}\mathcal{A}$. The canonical morphism
$\eta: T\rightarrow h^{sm}T$ is universal for morphisms of $T$
into smooth homotopy functors. It is evident that $h^{sm}T=T$ if
and only if $T$ is a smooth homotopy functor. The homotopization
$hT$ of a functor $T$ is defined similarly $A^{\infty (I)}$
replaced by $A^{I}$.

Since $K_n^{sm}$ are smooth homotopy functors, the homomorphism
$\alpha ^*_n(A):K_n(A)\rightarrow K_n^{sm}(A)$ yields a
homomorphism
\begin{equation}
h^{sm}\alpha ^*_n(A):h^{sm}K_n(A)\rightarrow
K_n^{sm}(A),\;\;\;\;n\geq 1,
\end{equation}
$A\in \mathrm{Ob}\mathcal{A}$.

\begin{thm}\label{smtha}
Let $A$ be a locally convex $k$-algebra. Then\\
(i) there is an isomorphism
$$
h^{sm}\alpha ^*_1(A): h^{sm}K_{1}(A)\xrightarrow{\approx}
K_1^{sm}(A),
$$
(ii) $\alpha ^*_1(A):K_1(A)\rightarrow K_1^{sm}(A)$ is an
isomorphism if $K_1(A)\rightarrow K_1(A^{\infty(I)})$ is an
isomorphism.
\end{thm}

\begin{proof}
(i) The commutative diagram
$$
\begin{array}{ccccccc}
 GL(\mathcal{I(A)}) & \rightarrow & GL(A) & \rightarrow & K_1^{sm}(A) & \rightarrow & 0 \\
  \downarrow &  & \downarrow &  & \parallel &  &  \\
  h^{sm}K_1(\mathcal{I}(A)) & \rightarrow & h^{sm}K_1(A) & \xrightarrow{h^{sm}\alpha ^*_1(A)} & K_1^{sm}(A) &
   \rightarrow & 0 \\
\end{array}
$$
with exact top row and vertical surjective homomorphisms implies
the exactness of the bottom row. The local convex $k$-algebra
$\mathcal{I}(A)$ is contractible, the trivial map
$0_{\mathcal{I}(A)}$ and the identity map $1_{\mathcal{I}(A)}$
being smooth homotopic with smooth homotopy $\delta
_A;\mathcal{I}(A)\rightarrow \mathcal{I}^2(A)$ between them. Since
$h^{sm}K_1$ is a smooth homotopy functor, one gets
$h^{sm}K_1(\mathcal{I}(A))=0$. Hence $h^{sm}\alpha ^*_1(A)$ is an
isomorphism.

(ii) Follows from (i), since in this case $h^{sm}K_{1}(A)=
K_{1}(A)$.
\end{proof}

It is easily checked that the smooth $K$-functors
$K^{sm}_n,\;\;n\geq1 $, are compatible with finite products of
locally convex $k$-algebras, namely there is a natural isomorphism
$$
K_n^{sm}\big(\prod_{n=1}^mA_i\big)\xrightarrow{\approx}\prod_{n=1}^mK_n^{sm}(A_i),\;\;\;
n\geq 1.
 $$

 \begin{prop}
 \label{smprpd}
 (i) Let $T:\mathcal{A}\rightarrow \mathcal{G}r$ be a functor such that
    $\mathrm{Im}(T(\varepsilon _0)\times T(\varepsilon _1))$ is a normal subgroup of
    $T(A)\times T(A),\;\;A\in \mathrm{Ob}\mathcal{A}$.
    Then there is an isomorphism
$$
\psi :h^{sm}T(A)\xrightarrow{\approx}\mathrm{Coker}(T(\varepsilon
_0)\times T(\varepsilon _1)),
$$
(ii) if an addition $T(0)=\{1\}$, then there is a surjection
$$
\mathrm{Coker}(T(\varepsilon _0\times \varepsilon _1))\rightarrow
h^{sm}T(A),\;\;\;A\in \mathrm{Ob}\mathcal{A}.
$$
 \end{prop}

\begin{proof}
Completely similar to the proof of Proposition 6 and Corollary 7
\cite{Inas1} and will be omitted.
\end{proof}

\begin{thm}
\label{smthb} Let $A$ be a locally convex $k$-algebra. Then for
fixed $i\geq 0$ there is an isomorphism
$$
\alpha _{i+1}:K_{i+1}(A)\rightarrow K^{sm}_{i+1}(A),
$$
if $K_{i+1}(A)\xrightarrow{\approx}K_{i+1}(A^{\infty (I)})$ and
$K_{j}(B^{\infty (I^{i-j})})\xrightarrow{\approx}K_{j}(B^{\infty
(I^{i-j+2})})$ with $B=A^{2^{l}}$ for all $0\leq j\leq i$, $0\leq
l\leq j-1$.
\end{thm}

\begin{proof}
We will give the proof for Swan's $K$-functors $K^s_*$. The case
$i=0$ is already proved (see Theorem \ref{smtha} (ii)).

Consider the short exact sequence of locally convex $K$-algebras
\begin{equation}
\label{smeqd} 0\rightarrow \Omega_{sm} (A)\rightarrow A^{\infty
(I)}\xrightarrow{\varepsilon} A\times A\rightarrow 0,
\end{equation}
$\varepsilon =\varepsilon _0\times \varepsilon _1$, which is a
$GL$-fibration with respect to the smooth path cotriple
$\mathbf{I}$. In effect, it is easy to see that the image of the
composites
$$
GL(\mathbf{I}_*(A))\rightarrow GL(\mathbf{I}_*(A^{\infty (I)}))
\rightarrow GL(\mathbf{I}_*(A))\times GL(\mathbf{I}_*(A)),
$$
$$
GL(\mathbf{I}_*(\mathbf{I}A))\rightarrow GL(\mathbf{I}_*(A^{\infty
(I)})) \rightarrow GL(\mathbf{I}_*(A))\times GL(\mathbf{I}_*(A)),
$$
is, respectively, the diagonal subgroup and the subgroup $1\times
GL(\mathbf{I}_*(A))$, where the map
$GL(\mathbf{I}_*(A))\rightarrow GL(\mathbf{I}_*(A^{\infty (I)}))$
is induced by the inclusion $i:A\rightarrow A^{\infty (I)}$ and
the map $GL(\mathbf{I}_*(A^{\infty (I)})) \rightarrow
GL(\mathbf{I}_*(A))\times GL(\mathbf{I}_*(A))$ is induced by
$\varepsilon _0\times \varepsilon _1$. Since the diagonal subgroup
and the subgroup $1\times GL(\mathbf{I}_*(A))$ generate the group
$GL(\mathbf{I}_*(A))\times GL(\mathbf{I}_*(A))$, it follows that
$GL(\mathbf{I}_*(\varepsilon _0\times \varepsilon _1))$ is
surjective.

Therefore the sequence \ref{smeqd} induces the long exact sequence
of smooth $K$-functors and yields the following exact sequence
$$
0\rightarrow K_i^{sm}(A^{\infty (I)})\rightarrow K_i^{sm}(A\times
A)\rightarrow K_{i-1}^{sm}(\Omega_{sm} (A))\rightarrow 0
$$
for $i\geq 1$.

 Now we will introduce auxiliary $K$-functors
$M_kK_i^s$, $k\geq 0$, $i\geq 1$, needed in the sequel.

Let us define the augmented simplicial rings $M_kF_*^+(A)$, $A\in
\mathrm{Ob} \mathcal{A}$, as follows:
$$
M_0F_*^+(A)=F_*^+(A),\;\;\;M_1F_*^+(A)=\mathrm{Ker}F_*^+(\varepsilon
_0\times \varepsilon _1).
$$
Any continuous $k$-homomorphism $f:A\rightarrow C$ of locally
convex $k$-algebras induces in a natural way a morphism
$M_1F_*^+(f):M_1F_*^+(A)\rightarrow M_1F_*^+(C)$ and one gets a
functor $M_1F_*^+$ from the category $\mathcal{A}$ to the category
of augmented simplicial rings. For $k>1$ we set
$M_kF_*^+(A)=M_1(M_{k-1}F_*^+(A))$, $A\in \mathrm{Ob}\mathcal{A}$.

Define $M_kK_i^s(A)=\pi _{i-2}GL(M_kF_*^+(A))$ for $i\geq 1$,
$k\geq 0$.

It is clear that any continuous $k$-homomorphism $f:A\rightarrow
C$ of locally convex $k$-algebras induces in a natural way a
morphism $M_kF_j^+(f):M_kF_j^+(A)\rightarrow M_kF_j^+(C)$ for all
$k\geq 0$ $j\geq 0$, where $\{M_kF_j^+(f)\}_j=M_kF_*^+(f)$. It is
also obvious that every augmented simplicial ring $M_kF_*^+(A)$,
$k\geq 0$, is a simplicial resolution of $\Omega_{sm}
^k(A)=\Omega_{sm} (\Omega_{sm} ^{k-1}(A))$, $\Omega_{sm} ^0(A)=A$,
in the category of rings.

Let $f:A\rightarrow C$ be a continuous $k$-homomorphism of locally
convex $k$-algebras such that there exists a $k$-linear continuous
linear map $f':C\rightarrow A$ with $ff'=1$. It is easily checked
that $f'$ induces $k$-linear continuous maps $\Omega_{sm}
(f'):\Omega_{sm} (C)\rightarrow \Omega_{sm} (A)$ and $f'^{\infty
(I)}:C^{\infty (I)}\rightarrow A^{\infty (I)}$ such that
$\Omega_{sm} (f)\Omega_{sm} (f')=1$ and $f^{\infty (I)}f'^{\infty
(I)}=1$. It is also clear that the $k$-linear continuous map $f'$
gives rise to a homomorphism of rings
$M_kF_j(f'):M_kF_j(C)\rightarrow M_kF_j(A)$ such that
$M_kF_j(f)M_kF_j(f')=\mathrm{identity}$ for all $k,j\geq 0$.
Therefore the induced homomorphism
$$
GL(M_kF_j(f)):GL(M_kF_j(A))\rightarrow GL(M_kF_j(C))
$$
is surjective for all $k,j\geq 0$ and one obtains a long exact
sequence
\begin{multline*}
...\rightarrow M_kK^s_{i+1}(C)\rightarrow
M_kK^s_{i}(A,\mathrm{Ker}(f))\rightarrow M_kK^s_{i}(A)\rightarrow
\\\rightarrow M_kK^s_{i} (C)\rightarrow
M_kK^s_{i-1}(A,\mathrm{Ker}(f))\rightarrow ...
\end{multline*}
for all $k\geq 0$, where
$$
M_kK^s_{i}(A,\mathrm{Ker}(f))=\pi
_{i-2}GL(\mathrm{Ker}(M_kF_*^+(f))),\;\;\;i\geq 1.
$$

We can apply this exact sequence to the homomorphism $\varepsilon
_0 \times \varepsilon _1:A^{\infty (I)}\rightarrow A\times A$ of
the sequence \ref{smeqd}. In effect, there is a $k$-linear
continuous map $\varepsilon ':A\times A\rightarrow A^{\infty (I)}$
given by $(a,a')\mapsto (1-t)a+ta'$ such that $(\varepsilon
_0\times \varepsilon _1)\varepsilon '=1_{A\times A}$. Thus the
sequence \ref{smeqd} yields the following long exact sequence
\begin{multline}
\label{smeqe} ...\rightarrow M_{j+1}K^s_{i+1}(A)\rightarrow
M_jK^s_{i+1}(A^{\infty (I)})\rightarrow M_jK^s_{i+1}(A\times
A)\rightarrow
\\\;\;\;\;\;\rightarrow M_{j+1}K^s_{i} (A)\;\;\rightarrow
\;\;M_jK^s_{i}(A^{\infty (I)})\;\;\rightarrow
\;\;M_jK^s_{i}(A\times
A)\rightarrow ...\\
...\rightarrow M_{j+1}K^s_{1} (A)\rightarrow M_jK^s_{1}(A^{\infty
(I)})\rightarrow M_jK^s_{1}(A\times A)
\end{multline}
for all $j\geq 0$.

The morphism $\alpha :\mathbf{F}\rightarrow \mathbf{I}$ of
cotriples induces a morphism $F_*(A)\rightarrow
\mathcal{I}_{*}(A)$ which yields a morphism $M_jF_*(A)\rightarrow
\mathcal{I}_{*}(\Omega _{sm}^j(A))$ for all $j\geq 0$ implying
isomorphisms
$$
h^{sm}(M_jK_1^{sm})(A)\approx K_1^{sm}(\Omega _{sm}^j(A))\approx
K_{j+1}^{sm}(A)
$$
for any locally convex $k$-algebra $A$ and $j\geq 0$.

By Proposition \ref{smprpd} one has a natural surjection
$$
\mathrm{Coker}(M_jK_i^{s}(\varepsilon _0\times \varepsilon
_1))\rightarrow h^{sm}(M_jK_i^{s}(A))
$$
for all $j\geq 0$, $i\geq 1$.

Thus for any $j$ satisfying $0\leq j\leq i-1$ the exact sequence
\ref{smeqe} provides the following commutative diagram:
\begin{equation}
\label{smeqf}
\begin{array}{ccc}
  M_{j+1}K_{i-j}^s(A)& \rightarrow & h^{sm}(M_jK_{i-j+1})(A) \\
  \downarrow &  & \downarrow \\
  K_{i-j}^{sm}(\Omega _{sm}^{j+1}(A)) & \xrightarrow{\approx} & K_{i+1}^{sm}(A) \\
\end{array}
\end{equation}
if $M_jK_{i-j}^{s}(A)\rightarrow M_jK_{i-j}^{s}(A^{\infty (I)})$
is an isomorphism, where the top homomorphism is surjective.

It follows that  if $M_jK_{i-j}^s(A)\rightarrow
M_jK_{i-j}^{s}(A^{\infty (I^2)})$ is an isomorphism for $0\leq
j\leq i-1$, then the diagram \ref{smeqf} holds also for $A^{\infty
(I)}$ and $0\leq j\leq i-1$. It is easy to verify that in this
case the diagram \ref{smeqf} yields the commutative diagram
\begin{equation}
\label{smeqg}
\begin{array}{ccc}
  h^{sm}(M_{j+1}K_{i-j}^s)(A)& \rightarrow & h^{sm}(M_jK_{i-j+1})(A) \\
  \downarrow &  & \downarrow \\
  K_{i-j}^{sm}(\Omega _{sm}^{j+1}(A)) & \xrightarrow{\approx} & K_{i+1}^{sm}(A) \\
\end{array}
\end{equation}
for all $j$ satisfying $0\leq j\leq i-1$, where the top
homomorphism is surjective.

Since $h^{sm}(M_iK_{1}^s)(A)\rightarrow K_{i+1}^{sm}(A)$ is an
isomorphism, the diagram
\begin{equation}
\label{smeqj}
\begin{array}{ccc}
  h^{sm}(M_{i}K_{1}^s)(A)& \rightarrow & h^{sm}(M_{i-1}K_{2}^s)(A) \\
  \downarrow &  & \downarrow \\
  K_{1}^{sm}(\Omega _{sm}^{i}(A)) & \xrightarrow{\approx} & K_{i+1}^{sm}(A) \\
\end{array}
\end{equation}
implies the isomorphism $h^{sm}(M_{i-1}K_{2}^s)(A)\rightarrow
K_{i+1}^{sm}(A)$.

Thus the diagram
\begin{equation}
\label{smeqi}
\begin{array}{ccc}
  h^{sm}(M_{i-1}K_{2}^s)(A)& \rightarrow & h^{sm}(M_{i-2}K_{3}^s)(A) \\
  \downarrow &  & \downarrow \\
  K_{2}^{sm}(\Omega _{sm}^{i-1}(A)) & \xrightarrow{\approx} & K_{i+1}^{sm}(A) \\
\end{array}
\end{equation}
implies the isomorphism $h^{sm}(M_{i-2}K_{3}^s)(A)\rightarrow
K_{i+1}^{sm}(A)$.

Continuing in this manner step by step we arrive to the diagram
\begin{equation}
\label{smeqh}
\begin{array}{ccc}
  h^{sm}(M_{1}K_{i}^s)(A)& \rightarrow & h^{sm}(K_{i+1}^s)(A) \\
  \downarrow &  & \downarrow \\
  K_{i}^{sm}(\Omega _{sm}(A)) & \xrightarrow{\approx} & K_{i+1}^{sm}(A) ,\\
\end{array}
\end{equation}
where the left vertical homomorphism is an isomorphism and the top
homomorphism is surjective. Hence $h^{sm}K^s_{i+1}(A)\rightarrow
K^{sm}_{i+1}(A)$ is an isomorphism, if $M_jK_{i-j}^s(A)\rightarrow
M_jK_{i-j}^{s}(A^{\infty (I^2)})$ is an isomorphism for $0\leq
j\leq i-1$.

Now conditions on the algebraic $K$-functors $K^s_j$ will be given
implying these conditions and therefore the isomorphism
$h^{sm}K^s_{i+1}(A)\rightarrow K^{sm}_{i+1}(A)$.

To this end let us consider the following commutative diagram with
exact rows
\begin{equation}
\label{smeqk}
\begin{array}{ccccccccc}
  &_{M_{j}K_{i+1}^{s}(A^{2})}& _{\rightarrow } & _{M_{j+1}K_{i}^{s}(A)} & _{\rightarrow }&
  _{M_{j}K_{i}^{s}
  (A^{\infty (I)})}& _{\rightarrow }& _{M_{j}K_{i}^{s}(A^2)}& \\
  &_{\downarrow }&  & _{\downarrow }&  & _{\downarrow }&  & _{\downarrow }&\\
 & _{M_{j}K_{i+1}^{s}((A^2)^{\infty (I^2)})} & _{\rightarrow }&
 _{M_{j+1}K_{i}^{s}(A^{\infty (I^2)})}
 & _{\rightarrow }& _{M_{j}K_{i}^{s}(A^{\infty (I^3)})}
 & _{\rightarrow }& _{M_{j}K_{i}^{s}((A^2)^{\infty (I^2)})}
 &\\
\end{array}
\end{equation}
induced by the exact sequence \ref{smeqe}.

The diagram \ref{smeqk} shows that for $j=0$ one has
$M_{1}K_{i}^{s}(A)\xrightarrow{\approx} M_{1}K_{i}^{s}(A^{\infty
(I^2)})$ if
\begin{enumerate}
    \item $K_{i+1}^{s}(A^{2})\xrightarrow{\approx}K_{i+1}^{s}((A^2)^{\infty
(I^2)})$,
    \item $K_{i}^{s}(A^{\infty (I)})\xrightarrow{\approx}K_{i}^{s}(A^{\infty
  (I^3)})$.
\end{enumerate}

Similarly the diagram \ref{smeqk} shows that for $j=1$ one has
$M_{2}K_{i}^{s}(A)\xrightarrow{\approx} M_{2}K_{i}^{s}(A^{\infty
(I^2)})$ if
$M_1K_{i+1}^{s}(A^{2})\xrightarrow{\approx}M_1K_{i+1}^{s}((A^2)^{\infty
(I^2)})$ and $M_1K_{i}^{s}
  (A^{\infty (I)})\xrightarrow{\approx}M_1K_{i}^{s}(A^{\infty
  (I^3)})$.
Thus one has an isomorphism
$M_{2}K_{i}^{s}(A)\xrightarrow{\approx} M_{2}K_{i}^{s}(A^{\infty
(I^2)})$ if
\begin{enumerate}
    \item $K_{i}^{s}(A^{\infty
(I^2)})\xrightarrow{\approx}K_{i}^{s}(A^{\infty
  (I^4)})$,
    \item $K_{i+1}^{s}
  ((A^2)^{\infty (I)})\xrightarrow{\approx}K_{i+1}^{s}((A^2)^{\infty
  (I^3)})$,
    \item $K_{i+2}^{s}
  (A^4)\xrightarrow{\approx}K_{i+2}^{s}((A^4)^{\infty
  (I^2)}).$
\end{enumerate}

By induction on $j$ it is easily checked that one has an
isomorphism
$$
M_jK^s_i(A)\xrightarrow{\approx}M_jK_i^s(A^{\infty
(I^2)}),\;\;j\geq 0,
$$
if the following isomorphisms hold:
$$
K_{i}^{s}
  (A^{\infty (I^j)})\xrightarrow{\approx}K_{i}^{s}(A^{\infty
  (I^{j+2})}),
  $$
  $$
 K_{i+1}^{s}
  ((A^2)^{\infty (I^{j-1})})\xrightarrow{\approx}K_{i+1}^{s}((A^2)^{\infty
  (I^{j+1})}),
  $$
  $$
     K_{i+2}^{s}
  ((A^4)^{\infty (I^{j-2})})\xrightarrow{\approx}K_{i+2}^{s}((A^4)^{\infty
  (I^j)}),
$$
.........................................................................\\
$$
   K_{i+j}^{s}
  (A^{2^j})\xrightarrow{\approx}K_{i+j}^{s}((A^{2^j})^{\infty
  (I^2)}),
$$
or shortly if
$$
K_{i+l}^{s}
  ((A^{2^l})^{\infty (I^{j-l})})\xrightarrow{\approx}K_{i+l}^{s}((A^{2^l})^{\infty
  (I^{j-l+2})})
$$
for all $l$ satisfying $0\leq l\leq j$.

We deduce that there is an isomorphism
$$
h^{sm}K^s_{i+1}(A)\rightarrow K_{i+1}^{sm}(A)
$$
if
$$
K_{j}^{s}
  ((A^{2^l})^{\infty (I^{i-j})})\xrightarrow{\approx}K_{j}^{s}((A^{2^l})^{\infty
  (I^{i-j+2})})
$$
for $1\leq j\leq i$, $0\leq l\leq j-1$.

To these conditions it suffices to add the isomorphism
$K_{i+1}^s(A)\rightarrow K^s_{i+1}(A^{\infty (I)})$ to obtain the
required isomorphism
$$
K_{i+1}^s(A)\rightarrow K^{sm}_{i+1}(A).
$$
This completes the proof.
\end{proof}

\begin{cor}
\label{smcra} Let $A$ be a locally convex $k$-algebra. Assume that
Quillen's $K$-functors $K_{i}$, $i\geq 1$, are compatible with
finite products of $A$ and $A^{\infty I}$ on themselves. Then
there is an isomorphism
$$
\alpha _*:K_*(A)\rightarrow K_*^{sm}(A),
$$
if $K_*(A)\xrightarrow{\approx}K_*(A^{\infty (I^n)})$ for all $n>
0$.
\end{cor}
\begin{thm}
\label{smthc} Let $A$ be a locally convex $k$-algebra. Then for
fixed $i\geq 0$ there is an isomorphism
$$
\beta
_{i+1}(A):K_{i+1}^{sm}(A)\xrightarrow{\approx}K_{i+1}^{top}(A),
$$
if
$$
 K_{i+1}^{sm}(A)\xrightarrow{\approx}K_{i+1}^{sm}(A^I),
 $$
 $$
    K_{j}^{sm}(A^{I^{i-j}})\xrightarrow{\approx}K_{j}^{sm}(A^{I^{i-j+2}})
$$
for all $0\leq j\leq i$.
\end{thm}
\begin{proof}
Similar to the proof of Theorem \ref{smthb}. In effect, the
sequence
\begin{equation}
\label{smeql} 0\rightarrow \Omega A\rightarrow
A^I\xrightarrow{\overline{\varepsilon}_0\times
\overline{\varepsilon}_1}A\times A\rightarrow 0
\end{equation}
is a $GL$-fibration with respect to the smooth path cotriple that
follows from the natural commutative diagram
$$
\begin{array}{ccc}
  A^{\infty (I)} & \xrightarrow{\varepsilon _0\times
\varepsilon _1} & A\times A \\
  \sigma \downarrow \;\;\;\;&  & \parallel \\
  A^{I} & \xrightarrow{\overline{\varepsilon}_0\times
\overline{\varepsilon}_1}  & A\times A ,\\
\end{array}
$$
where $\Omega (A)= Ker \tau_{A}, \tau_{A}: J(A)\rightarrow A $ and
$\sigma $ is the natural inclusion, the continuous
$k$-homomorphism $\varepsilon _0\times \varepsilon _1$ being a
$GL$-fibration with respect to the smooth path cotriple
$\mathbf{I}$. Replace the auxiliary $K$-functors $M_jK_i^s(-)$ by
the $K$-functors $K_{i}^{sm}(\Omega^{j}_{sm}(-))$, use the natural
isomorphisms
$$
K^{sm}_{*}((A\times A)^{I^{n}})\approx K^{sm}_{*}(A^{I^{n}})\times
K^{sm}_{*}(A^{I^{n}}), n\geq 0 ,
$$
the canonical isomorphism $hK^{sm}_{1}\approx K^{top}_{1}$, where
$hK^{sm}_{1}$ is the continuous homotopyzation of the functor
$K^{sm}_{1}$, and note that in this case
$M_{j}\mathcal{I}^{+}_{*}(A)=
Ker\;M_{j-1}\mathcal{I}^{+}_{*}({\overline\varepsilon}_{0}\times
{\overline\varepsilon}_{1})$ is the $\mathcal{I}$-projective
cotriple resolution of $\Omega^{j}_{sm}(A)$ for all $j\geq 1$.
\end{proof}

A Fr\'echet $k$-algebra whose seminorms possess the multiplicative
property will be called Fr\'echet-Michael $k$-algebra. It is easy
to see that the smooth and continuous cotriples are both defined
on the category of Fr\'echet-Michael $k$-algebras and their
continuous $k$-homomorphisms.

\begin{thm} Let $A$ be a complex Fr\'echet-Michael algebra. Then
there is an isomorphism
$$
\beta_{n}(A): K_{n}^{sm}(A)\xrightarrow{\approx}K_{n}^{top}(A)
$$
 for all $n\geq 0$.
 \end{thm}
 \begin{proof} The case $n=0$ is trivial and it suffices to prove the
 theorem for the unital case.

  First it will be shown that
 $\beta_{1}(A): K_{1}^{sm}(A)\xrightarrow{\approx} K_{1}^{top}(A)$.

Let $f: A\rightarrow B$ be a surjective homomorphism of arbitrary
complex Fr\'echet-Michael algebras. According to Lemma 1.14 [15]
the induced map $GL^{0}(f): GL^{0}(A) \rightarrow GL^{0}(B)$ is
surjective, where $GL^{0}(-)$ denotes the component of the unit of
$GL(-)$. Therefore the canonical map $\tau_{A}:
\mathcal{I}(A)\rightarrow A$ induces a surjective map
$GL^{(0)}(\tau_{A}): GL^{(0)}(\mathcal{I}(A))\rightarrow
GL^{(0)}(A)$. On the other hand $\mathcal{I}(A)$ is continuously
contractible. In effect, take the map $\bar{\delta}_{A}:
\mathcal{I}(A) \rightarrow (\mathcal{I}(A)^{\infty(I)}$ which is
the composite of
$$
\mathcal{I}(A)\xrightarrow{{\delta}_{A}}\mathcal{I}^{2}(A)\rightarrow(\mathcal{I}(A)^{\infty(I)}.
$$
One has the following commutative diagram
\begin{equation}
\begin{array}{ccccc}
   \mathcal{I}(A)& \xrightarrow{\bar{\delta}_A} & \mathcal{I}(A)^{\infty (I)} & \xrightarrow{\varepsilon _i} & \mathcal{I}(A) \\
  \| &  & \theta \downarrow \;\;&  & \|\\
  \mathcal{I}(A)& \xrightarrow{\theta \bar{\delta}_A}& \mathcal{I}(A)^I& \xrightarrow{\bar{\varepsilon}_i} & \mathcal{I}(A), \\
\end{array}
\end{equation}
$i=0,1$, where $\bar{\varepsilon}_{i}= f(i)$ and $\theta$ is the
natural inclusion. Diagram (1.18) implies the equalities
$\bar{\varepsilon}_{0}\theta \bar{\delta}_{A}= 0$ and
$\bar{\varepsilon}_{1}\theta \bar{\delta}_{A}= 1_{\mathcal{I}(A)}$
showing that $\mathcal{I}(A)$ is a contractible Fr\'echet-Michael
algebra. Thus, $GL(\mathcal{I}(A))$ is a connected space and
therefore one has a surjection $GL(\tau_{A}):
GL(\mathcal{I}(A))\rightarrow GL^{(0)}(A)$. We deduce
$$
K_{1}^{sm}(A)= GL(A)/ImGL(\tau_{A})= GL(A)/GL^{(0)}(A)=
K_{1}^{top}(A)
$$
for any complex Fr\'echet-Michael algebra $A$.

To obtain the required isomorphism in higher dimensions it should
be noted that any surjective homomorphism $f:A\rightarrow B$ of
complex Fr\'echet-Michael algebras is a $GL$-fibration with
respect to the continuous path cotriple. In fact by Lemma 1.9 [15]
one has surjections $J_{n}(A)\rightarrow J_{n}(B)$, $n\geq 0$,
induced by $f$, and therefore the induced map $GL^{(0)}(J_{n}(f)):
GL^{(0)}(J_{n}(A))\rightarrow GL^{(0)}(J_{n}(B))$ is surjective
for all $n\geq 0$. Since $J_{n}(A)$ and $J_{n}(B)$ are
contractible complex Fr\'echet-Michael algebras, $GL(J_{n}(A))$
and $GL(J_{n}(B))$ are connected spaces. We conclude that the
homomorphism $GL(J_{n}(f))$ is surjective for all $n\geq 0$ and by
definition $f$ is a $GL$-fibration with respect to the path
cotriple $\textbf{J}$.

The smooth and continuous contractibility of $\mathcal{I}(A)$
implies $K_{n}^{sm}(\mathcal{I}(A))=
K_{n}^{top}(\mathcal{I}(A))=0$ for $n> 0$ and any complex
Fr\'echet-Michael algebra $A$. It remains to apply to the short
exact sequence (1.6) of complex Fr\'echet-Michael algebras the
long exact sequence of smooth and topological $K$-functors
respectively and one obtains finally the following isomorphisms
$$
K_{n}^{sm}(A)\approx K_{1}^{sm}(\Omega_{sm}^{n-1}(A))\approx
K_{1}^{top}(\Omega_{sm}^{n-1}(A))\approx K_{n}^{top}(A)
$$
for $n> 1$ and any complex Fr\'echet-Michael algebra $A$. This
completes the proof.
\end{proof}
\begin{cor} Let $A$ be a complex Fr\'echet-Michael algebra such
that $GL_{1}(A^{+})$ is open in $A^{+}$. Then the $K$-groups
$K_{n}^{sm}(A)$, $K_{n}^{top}(A)$ and $RK_{n}(A)$ (defined in
[15]) are isomorphic for all $n\geq 0$. In particular, on the
category of complex Banach algebras the smooth $K$-theory is
isomorphic to the well-known topological $K$-theory of\underline{}
Banach algebras.
\end{cor}

We can interpret smooth $K$-functors $K^{sm}_{*}$ as "smooth
homotopy groups" of the general linear group.

The group $GL_{m}(A)= \mathrm{Ker}(GL_{m}(A^{+})\rightarrow
GL_{m}(k))$ with topology induced from $M_{m}(A^{+})$ in general
is not a topological group, since the map
$(-)^{-1}:GL_{m}(A)\rightarrow GL_{m}(A)$ is not always
continuous, depending on the locally convex $k$-algebra $A$.

Denote by $R^{n}GL_{m}(A)$ the group of all continuous maps
$f:I^{n}\rightarrow GL_{m}(A)$ satisfying the following
conditions:

(i)  $f(t_{1},t_{2},...,t_{n})=1$ if any $t_{i}=0$;

(ii) the composite maps $\varphi_{ij}.f:I^{n}\rightarrow A$ are
smooth maps for all $1\leq i,j\leq m$, where the map
$\varphi_{ij}:GL_{m}\rightarrow A$ is given by
$\varphi_{ij}(M)=a_{ij}$ with $a_{ij}$ the element of $A$ staying
at the intersection of the i-th row and the j-th column of the
matrix $M$;

(iii) the map $f^{-1}:I^{n}\rightarrow GL_{m}(A)$ given by
$f^{-1}(t_{1},t_{2},...,t_{n})= (f(t_{1},t_{2},...,t_{n}))^{-1}$
is continuous and satisfies condition (ii).

Using the homotopy groups of the simplicial group
$R_{*}GL_{m}(A)$, whose boundary and degeneracy maps coincide with
those of $GL_{m}(I_{*}(A))$, we define by setting
$$
\pi^{sm}_{n}(GL_{m}(A))= \pi_{n-1}R_{*}GL_{m}(A),\;n\geq 2,
$$
and by the exact sequence
$$
0\rightarrow \pi^{sm}_{1}GL_{m}(A)\rightarrow
\pi_{0}R_{*}GL_{m}(A)\rightarrow GL_{m}(A)\rightarrow
\pi^{sm}_{0}GL_{m}(A)\rightarrow 0.
$$

Finally, define
$$
\pi^{sm}_{n}GL(A)= \underrightarrow{\lim}_m \pi^{sm}_{n}GL_{m}(A)
$$
for $n\geq 0$. It is obvious that one has
$$
K^{sm}_{n}(A)=\pi^{sm}_{n-1}GL(A),\; n\geq 1.
$$

\section{The Cohen-Hewitt factorization for Fr\'echet $k$-algebras}

 It is intended to extend the Cohen-Hewitt theorem factorization theorem to the
category of Fr\'echet $k$-algebras which will be used in the
sequel.

We will assume that $A$ is a $k$-algebra, $L$ is a left $A$-module
equipped with a $k$-linear space structure and
$$
(t\mu )x=t(\mu x)=\mu (tx)
$$
for $t\in k$, $\mu \in A$ and $x\in L$.

Particularly we will consider $A$ as a Banach $k$-algebra and $L$
as a Fr\'echet space. Since $L$ is a Fr\'echet space, its topology
is given by an increasing sequence of seminorms ${||\cdot ||_n}$
satisfying the following conditions:
\begin{itemize}
    \item For any seminorm $||\cdot ||_n$ there exist a number $C\geq
    1$ independent of $n$ and a seminorm $||\cdot||_{m(n)}$ such that
    \begin{equation}
    \label{chweqa}
||\mu x||_n\leq C||\mu||\cdot||x||_{m(n)}
    \end{equation}
\end{itemize}
for all $\mu\in A$, $x\in L$ and $n$. It is obvious that the
module structure $A\times L\rightarrow L$ is jointly continuous.
Therefore $L$ is a topological $A$-module which will be called a
Fr\'echet module over a Banach $k$-algebra.

Thus a metric $\rho$ compatible with the Fr\'echet topology can be
defined on $L$ given by
\begin{equation}
\label{chweqb} \rho (x,y)=\sum
_{n=1}^{\infty}\frac{1}{2^n}\frac{||x-y||_n}{1+||x-y||_n}.
\end{equation}

\begin{defn}
Let $L$ be a Fr\'echet module over a Banach $k$-algebra $A$. It
will be said that $L$ possesses a left bounded approximate unit
(bounded by a positive constant $d$) in $A$, if for any finite
subset $\{\mu _1,...,\mu _m\}\subset A$, any element $x\in L$ and
$\varepsilon
>0$ there exists $\nu \in A$ such that
\begin{equation}
\label{chweqc} ||\nu||\leq d,\;\;\;\;\;||\nu \mu _i-\mu
_i||<\varepsilon,\;\;\;\;\; \rho (\nu x;x)<\varepsilon .
\end{equation}
\end{defn}

This is equivalent to the existence of a bounded direct set $\{\nu
_{\alpha }\}$ in $A$ such that
$$
\underrightarrow{\lim}_\alpha (\nu _\alpha \mu) =\mu \;\;\;\;\;\;
\text{and}\;\;\;\;\;\;\underrightarrow{\lim}_\alpha (\nu _\alpha
x)=x,
$$
$\mu \in A,\; x\in L$.

Denote by $A^+$ the Banach algebra $k$-algebra $A+k=\{(a,t)|\;a\in
A,\;t\in k\}$ with usual sum, with product given by
$(a,t)(a',t')=(aa'+t'a+ta',tt')$ and with norm
$||(a,t)||=||a||+|t|$.

For $\mu \in A$ with $||\mu||\leq d$ and $d\geq 1$, define an
element $\varphi (\mu)$ of $A^+$ by
\begin{equation}
\label{chweqd} \varphi (\mu)=\frac{2d+1}{2d} \bigg(1+\sum
_{k=1}^{\infty}(-1)^k(2d)^{-k}\mu ^k\bigg).
\end{equation}

Then one has (see \cite{Hewi})
\begin{equation}
\label{rtr} \varphi (\mu )=\bigg[\frac{2d}{2d+1}+\frac{1}{2d+1}\mu
\bigg]^{-1}\;\;\text{and}\;\;\;\frac{2}{3}+d^{-1}\leq ||\varphi
(\mu)||\leq 2+d^{-1}.
\end{equation}

The following lemma  is a generalization of Lemma 2.1 \cite{Hewi}.
\begin{lem}\label{chwla}
Let $L$ be a Fr\'echet module over a Banach $k$-algebra. Then for
any element $\mu \in A$ with $||\mu||\leq d\geq 1$ one has
\begin{equation}
\rho (\varphi (\mu)x;x)\leq C(2+d^{-1})\rho (\mu x;x).
\end{equation}
\end{lem}
\begin{proof}
 The following inequality is a consequence of the inequality
(\ref{chweqa}) and the properties (\ref{rtr}):
\begin{multline*}
\label{strs}
||\varphi(\mu )x-x||_n =\\
=||\varphi(\mu )x-\varphi(\mu )\cdot
\bigg(\frac{2d}{2d+1}+\frac{1}{2d+1}\mu \bigg) x||_n\leq
\\\leq C(2d+1)^{-1}||\varphi(\mu )||\cdot||\mu x-x||_{m(n)}\leq \\
\leq C(2+d^{-1})||\mu x-x||_{m(n)}.
\end{multline*}
Since $(2+d^{-1})C\geq 1$, one has
$$
\frac{||\varphi (\mu )x-x||_n}{1+||\varphi (\mu )x-x||_n}\leq
C(2+d^{-1})\cdot \frac{||\mu x-x||_{m(n)}}{1+||\mu x-x||_{m(n)}}.
$$
Therefore
$$ \sum
_{n=1}^{\infty}\frac{1}{2^n}\cdot \frac{||\varphi (\mu
)x-x||_n}{1+||\varphi (\mu )x-x||_n}\leq C(2+d^{-1})\cdot \sum
_{i=1}^{\infty}\frac{1}{2^i}\cdot \frac{||\mu x-x||_i}{1+||\mu
x-x||_i}.
$$
Recall the definition of the metric $\rho $ to obtain the required
inequality.
\end{proof}

Let $d >0$ and $\{\mu _n\}$ be a sequence of elements of $A$ such
that $||\mu _n||\leq d,\;\;\;n=1,2,...$. Define an element $\sigma
_n \in A^+$ by
\begin{equation}\label{chweqe}
\sigma _n=\sum_{k=1}^{n}(2d)^{k-1}(2d+1)^{-k}\mu_k+(2d)^n(2d+1)^n
\end{equation}
for $n=1,2,...,$ and $\sigma _0=1$.

It is known that $\sigma _n$ is invertible in $A^+$ \cite{Hewi}.

\begin{prop}\label{chwpa}
Let $L$ be a Fr\'echet module over a Banach $k$-algebra $A$ having
a bounded approximate unit in $A$, and let $z\in L$. Let $d$ be an
upper bound on the norms of approximate units in $A$. Then for
$\varepsilon >0$ there exists a sequence $\mu _n$ of elements of
$A$ with $||\mu _n||\leq d$, such that the following inequality
holds
\begin{equation}\label{chweqf}
\rho(\sigma _n^{-1}\cdot z;\sigma _{n-1}^{-1}\cdot z)\leq
\frac{\varepsilon}{2^n}.
\end{equation}
\end{prop}

\begin{proof}
Since $A$ has an approximate unit bounded by $d>0$, for
$\varepsilon >0$ there exists $\mu _1 \in A$ with $||\mu _1||$
such that
$$
\rho (\mu _1z;z)\leq(2C(2+d^{-1}))^{-1}\varepsilon.
$$
Then
$$
\sigma _1=\frac{1}{2d+1}\mu _1+\frac{2d}{2d+1}
$$
and $\sigma _1^{-1}=\varphi (\mu _1)$.

By Lemma \ref{chwla}
$$
\rho (\sigma _1^{-1}z;\sigma _0^{-1})\leq C(2+d^{-1})\rho (\mu
_1z;z)\leq \frac{\varepsilon}{2}.
$$
That means Proposition \ref{chwpa} is true for $n=1$.

Suppose Proposition \ref{chwpa} holds for $n=m$ and it will be
shown that it remains true for $n=m+1$.

Let us consider an element $\mu '\in A$ such that $||\mu '||\leq
d$ and define for $\varepsilon >0$
\begin{multline}\label{chweqg}
\sigma _{m+1}'=\Sigma_{k=1}^m(2d)^{k-1}(2d+1)^{-k}\mu
_k+(2d)^{m}(2d+1)^{-m-1}\mu '+\\+(2d)^{m+1}(2d+1)^{-m-1}.
\end{multline}
It is clear that \ref{chweqg} yields $\sigma _{m+1}$ when $\mu '$
is replaced by $\mu _{m+1}$. Rewrite \ref{chweqg} in the following
form
\begin{equation}\label{chweqj}
\sigma _{m+1}'=(\frac{1}{2d+1}\mu '+\frac{2d}{2d+1})\tau _m,
\end{equation}
where
$$
\tau _m=\sum _{k=1}^m(2d)^{k-1}(2d+1)^{-k}\varphi(\mu ')\mu
_k+(2d)^m(2d+1)^{-m}.
$$
The element $\mu '$ can be chosen such that $\tau _m$ becomes
invertible in $A^+$. In effect, since the group of invertible
elements in $A^+$ is an open subset and the map $x\mapsto x^{-1}$
is an homeomorphism, for $||\mu '\mu _k-\mu _k||$ sufficiently
small (this can be realized $A$ having a bounded approximate unit)
$||\tau _m-\sigma _m||$ is sufficiently small and therefore $\tau
_m$ is invertible. It follows that $||\tau _m^{-1}-\sigma
_m^{-1}||$ is also arbitrary small. Further, $\mu '$ can be chosen
such that $\rho (\mu 'z;z)$ will be sufficiently small. For such
$\mu '$ one has
\begin{multline}
\label{chweqi} \rho((\sigma _{m+1}')^{-1}z;\sigma
_m^{-1}z)=\rho(\tau _{m}^{-1}\varphi (\mu ')z;\sigma _m^{-1}z)
\leq \\\rho(\tau _{m}^{-1}\varphi (\mu ')z;\tau
_m^{-1}z)+\rho(\tau _{m}^{-1}z;\sigma _m^{-1}z).
\end{multline}
By Lemma \ref{chwla} combined with above, $\rho (\varphi (\mu
')z;z)$ is arbitrary small. Hence, since $\tau
_m^{-1}:L\rightarrow L$ is a homeomorphism given by $x\mapsto \tau
_m^{-1}x$, this implies that
$$
\rho((\tau _{m}^{-1}\varphi (\mu ')z;\tau _m^{-1}z)
$$
is also arbitrary small.

We have also to show that the second summand $\rho((\tau
_{m}^{-1}z;\sigma _m^{-1}z)$ of \ref{chweqi} is arbitrary small.
In fact, any $z\in L $ induces a continuous map $A\rightarrow L$
of metric spaces given by $a\mapsto az$. Hence, if $||\tau
_{m}^{-1}z-\sigma _m^{-1}z||$ is sufficiently small , then
$\rho((\tau _{m}^{-1}z;\sigma _m^{-1}z)$ is also sufficiently
small. Therefore the element $\mu _{m+1}=\mu '$ can be  chosen
such that
\begin{equation}\label{chweqk}
\rho(\sigma _{m+1}^{-1}\cdot z;\sigma _{m}^{-1}\cdot z)\leq
\frac{\varepsilon}{2^{m+1}}.
\end{equation}
This completes the proof.
\end{proof}

Now we are ready to prove the theorem generalizing Cohen-Hewitt's
result (Theorem 2.5 \cite{Hewi})

\begin{thm}\label{chwta}
Let $L$ be a Fr\'echet module over a Banach $k$-algebra $A$ having
a bounded approximate unit bounded by a number $d>0$. Then for
$z\in L$ and $\varepsilon >0$ there exist elements $\sigma \in A$
and $y\in L$ with the following properties

(i) $z=\sigma y$;

(ii)  $y\in \overline{A\cdot z}$\;\;\;(closure in L);

(iii) $\rho (y;z)\leq \varepsilon $;

(iv)  $||\sigma ||\leq d$.
\end{thm}

\begin{proof}
Take $\sigma _n$ defined as above (see equality \ref{chweqe})
which is invertible in $A^+$. Consider the elements $y_n=\sigma
_nz,\;\;n=1,2,...$ By summing the inequality
\begin{equation}\label{chweqh}
\rho (\sigma _n^{-1}z;\sigma _{n-1}^{-1}z)\leq
\frac{\varepsilon}{2^n}
\end{equation}
from $m$ to $m+k$, one obtains the inequality
$$
\rho (y_m;y_{m+k})\leq \frac{\varepsilon}{2^m},
$$
showing that the sequence $\{y_n\}$ is a Cauchy sequence and
denote $y=\underrightarrow{\lim}y_n$.

Since $L$ possesses an approximate unit, one has $z\in
\overline{Az}$ implying that the elements $y_n$ and $y$ belong to
$\overline{Az}$. By summing the inequality \ref{chweqh} from $1$
to $m$, and taking into account that $\sigma _0=1$, one gets the
inequality $ \rho (y_m;z)<\varepsilon$ and therefore
$$
\rho (y;z)\leq\varepsilon .
$$
By the definition of $\sigma _n$ it follows that
$\underrightarrow{\lim}\sigma_n$ exists in $A^+$ and in fact it is
an element $\sigma$ of $A$,
$$
\sigma =\Sigma _{k=1}^{\infty}(2d)^{k-1}(2d+1)^k\mu _k.
$$
Clearly $||\sigma||\leq d$. Since the module map is jointly
continuous, this implies
$$
z=\underrightarrow{\lim}(\sigma _n\cdot
y_n)=\underrightarrow{\lim}\sigma _n \cdot\underrightarrow{\lim}
y_n=\sigma \cdot y.
$$
This completes the proof.
\end{proof}

\section{Homotopy invariance in $K$-homology and Higson's Theorem \label{homin}}

The purpose of this section is, according to the homotopy
invariance of $K$-homology, to present Higson's homotopy
invariance theorem for both real and complex cases. Higson's
theorem asserts
\begin{thm}
\label{higs} (\cite{Higs} (Theorem 3.2.2)) Let $E$ be a stable and
split additive functor from an admissible sub-category
$\mathcal{S}$ of the category of complex $C^{*}$-algebras into the
category of abelian groups. Then it is homotopy invariant.
\end{thm}

 Higson's proof is a consequence of the following proposition
(cf. Theorem 3.1.4 in \cite{Higs}.)
\begin{prop}
\label{frp} Let $E$ be a functor from the category $\mathcal{S}$
into the category $Ab$ admitting a pairing with the set of
Fredholm $B$-pairs, $B\in Ob\mathcal{S}$. Then $E$ is a homotopy
functor.
\end{prop}

We recall the definition of a pairing of a functor $E:S\rightarrow
Ab$ with the set of Fredholm pairs, defined in \cite{Higs}, where
$\mathcal{S}$ is an admissible subcategory of the category
$C^*$-algebras. A subcategory $\mathcal{S}$ of the category of
$C^*$-algebras and $*$ -homomorphisms is said to be admissible if
\begin{enumerate}
    \item $k$ belongs to $\mathcal{S}$;
    \item if $A$ belongs to $\mathcal{S}$, then so is $A\otimes k^{I}$;
    \item if $A$ belongs to $\mathcal{S}$, then so is $A\otimes
\mathcal{K}$;
    \item if $0\rightarrow A\rightarrow B\rightarrow
C\rightarrow 0$ is a split exact sequence with $A$ and $C$ in
$\mathcal{S}$, then so is $B$.
\end{enumerate}
Here $"\otimes " $ is the minimal $C^*$-tensor product. A Fredholm
$B$-pair is a pair $(\varphi,\psi)$ of $*$-homomorphisms from $B$
into $\mathcal{L}_k(\mathcal{H})$ such that $\varphi (b)-\psi
(b)\in \mathcal{K}(\mathcal{H})$ for any $b\in B$, where
$\mathcal{H}$ is a countably generated Hilbert space over $k$,
Here $\mathcal{K}(\mathcal{H})$ is the $C^{*}$-algebra of compact
operators. A pairing of $E$ with the set of Fredholm $B$-pairs
assigns to each Fredholm $B$-pair $(\varphi ,\psi )$ a
homomorphism $\times (\varphi ,\psi ):E(A\otimes B)\rightarrow
E(A\otimes k)$ for any $A,B\in Ob\mathcal{S}$ with the following
properties:

\begin{enumerate}
\item  \underline{Functoriality. } If $(\varphi ,\psi )$\ is a
Fredholm $B'$-pair and if $f:B\rightarrow B'$ is a
$*$-homomorphism of $\mathcal{S}$, then the diagram
$$
\begin{array}{ccc}
E(A\otimes B) & \xrightarrow{\times (\varphi f,\psi f)}&
E(A\otimes k)\\
_{E(\mathrm{id}_A\otimes f)}\downarrow \;\;\;& &\parallel \\
E(A\otimes B') & \xrightarrow{\times (\varphi ,\psi )}& E(A\otimes
k)
\end{array}
$$
commutes.

\item  \underline{Additivity.} If $(\varphi ,\chi )$\ and $(\chi
,\psi )$\ are Fredholm $B$-pairs, then
$$
\times (\varphi ,\chi )+\times (\chi ,\psi )=\times (\varphi ,\psi
).
$$

\item  \underline{Stability.} If $(\varphi ,\psi )$\ is a Fredholm
$B$-pair and $\eta :B\rightarrow \mathcal{L}_k(\mathcal{H})$ is
any $*$-homomorphism, then
$$
\times (\varphi ,\psi )=\times \left( \left(
\begin{array}{cc}
\varphi  & 0 \\
0 & \eta
\end{array}
\right) ,\left(
\begin{array}{cc}
\psi  & 0 \\
0 & \eta
\end{array}
\right) \right) .
$$

\item  \underline{Non-degeneracy.} If $(e,\theta)$\ is a Fredholm
$B$-pair, $e:k\rightarrow \mathcal{K}(\mathcal{H})$ maps $1\in k$
to $p$, where $p$ is a rank one projection in
$\mathcal{K}(\mathcal{H})$ and $\theta $ is the zero homomorphism,
then
$$
\times (e,\theta):E(A\otimes k)\rightarrow E(A \otimes k)
$$
is the identity morphism.

\item  \underline{Unitary equivalence.} If $U\in
\mathcal{L}(\mathcal{H})$ is a unitary operator, then
$$
\times (\varphi ,\psi )=\times (U\varphi U^{*},U\psi U^{*}).
$$

\item  \underline{Compact perturbations.} If $U\in
\mathcal{L}(\mathcal{H})$ is a unitary operator equal to the
identity modulo compacts, then
$$
\times (\varphi ,U\varphi U^{*})=0.
$$
\end{enumerate}

Let $\hom(E(A\otimes -),E(A\otimes k))$ be a represented
contravariant functor from an admissible category to the category
of abelian groups.

The following proposition is crucial for the proof of proposition
\ref{frp}.

\begin{prop}
\label{higpr} Let $E$ be a functor from an admissible category of
$C^*$-algebras to the category of abelian groups and assume there
is a pairing of $E$ with the set of Fredholm $B$-pairs. Then there
is a natural transformation of functors
\begin{equation}
\phi :KK( -,k)\rightarrow \hom(E(A\otimes -),E(A\otimes k)).
\end{equation}
sending the identity of $KK(k,k)$ to the identity of
$\hom(E(A\otimes k),E(A\otimes k))$.
\end{prop}

Before showing this proposition, which will be based on an
investigation in $KK$-theory due to J. Cuntz and G. Scandalis
\cite{CuSk}, we need some remarks about functional calculus for
both real and complex $C^*$-algebra cases. Since we have not seen
functional calculus for real $C*$-algebras in the literature, we
will explain what we mean. The functional calculus for a
self-adjoint element $x$ in a complex $C^{*}$-algebra $A$ is the
$*$-monomorphism $\Phi :C(\mathrm{sp}x)\rightarrow A$, defined by
$id_{\mathrm{sp}x}\mapsto x$. Let $A$ be a real $C^{*}$-algebra
and consider the complex involutive algebra $A\otimes _RC$ with
involution $(a\otimes c)^{*}=a^{*}\otimes \bar c$. Then there
exist a $C^{*}$-norm on $A\otimes _RC$ and a canonical
$*$-embedding $A\hookrightarrow A\otimes _RC$ defined by $a\mapsto
a\otimes 1$ (cf. Theorem 2 and Corollary 2 in \cite{Palm}. Let
$r\in A$ be a self-adjoint element in the real $C^*$-algebra $A$
and $R(sp(r\otimes 1))$ be the real algebra of continuous real
functions on $sp(r\otimes 1)$, then the map $id_{sp(r\otimes
1)}\mapsto r$ defines a homomorphism $\psi:R(sp(r\otimes
1))\rightarrow A$ such that diagram
$$
\begin{array}{ccc}
  R(sp(r\otimes 1)) & \xrightarrow{\psi} & A \\
  \bigcap &  & \bigcap \\
  C(sp(r\otimes 1)) & \xrightarrow{\Psi} & A\otimes _RC\\
\end{array}
$$
commutes. This implies that $\psi:R(sp(r\otimes 1))\rightarrow A$
is a monomorphism which is called real functional calculus of the
self-adjoint element $r$ in a real algebra $A$.

 Now we return to the proof of Proposition \ref{higpr}.
\begin{proof}
Since the functional calculus exists for real $C^*$-algebras, the
technique of 17.4 and 17.6 in \cite{Blac} can be applied not only
for complex $C^*$-algebras, but also for real $C^*$-algebras.
Further, $KK(-,k)$ can be replaced naturally by $KK_c(-,k)$
(\cite{Blac}, Theorem 17.10.7). Now we construct a natural
transformation
\begin{equation}
\vartheta :KK_c( -,k)\rightarrow \hom(E(A\otimes -),E(A\otimes
k)).
\end{equation}
as follows. Recall  that if $(E,\varphi ,F)$ is a Kasparov
$(B,k)$-module, according to the results of subsections 17.4 and
17.6 in \cite{Blac}, one can construct a Fredholm pair $(\varphi
_0,\varphi _1)$ having the following properties (see Chapter 17,
6.2-6.3 in \cite{Blac}):
\begin{enumerate}
    \item if $(E,\varphi ,F)$ is degenerated, then $\varphi _0=\varphi
    _1$;
    \item Unitary equivalence corresponds to conjugation between $\varphi _0$ and $\varphi
    _1$ by the same unitary operator.
    \item "compact perturbation" corresponds to conjugation of $\varphi
    _1$ by a unitary operator which is a compact perturbation of the
    identity.
\end{enumerate}
Comparing properties (3), (5), (6) of the pairing with the
properties (1),(2),(3), one immediately concludes that there is a
map $\theta _B:KK_c(B,k)\rightarrow \hom(E(A\otimes B),E(A\otimes
k))$ defined by
$$
(E,\varphi ,F)\mapsto \times (\varphi _0,\varphi _1).
$$
Properties (1) and (2) of the pairing guarantee that $\{\theta
_B\}$ is a natural transformation of functors with values in the
category of abelian groups. The last requirement of the
proposition follows from the property (4) of the pairing.
\end{proof}

Now we are ready to prove Proposition \ref{higpr}.

\begin{proof} Since $K$-homology has the homotopy invariance property
with respect to the first variable, the diagram
\begin{equation}
\begin{array}{ccc}
  KK(k,k) & \xrightarrow{e_0=e_1} & KK(k[0,1],k) \\
  _{\theta _k}\downarrow &  & \downarrow _{\theta _{k[0,1]}}\\
  \hom(E(A\otimes k),E(A,\otimes k) & \xrightarrow{\bar{e}_0\sim \bar{e}_1} & \hom(E(A\otimes
  k[0,1]),E(A\otimes k ) \\
\end{array}
\end{equation}
commutes, where $e_0$, $e_1$, $\bar{e}_0$, $\bar{e}_1$ are induced
by the evolution maps $ev_0,\;ev_1:k[0,1]\rightarrow k$ at $0$ and
$1$ respectively. Let $\iota$ be a class in $KK(k,k)$ of the
Fredholm pair $(e,\theta)$ having property (4) of the Fredholm
pairing. Then $\theta _k(\iota)=id_{E(A)}$. Since $\theta
_{k[0,1]} e_0=\theta _{k[0,1]} e_1$, this implies
$$
\epsilon _0=\bar{e}_0(id_E(A))=\bar{e}_0(\theta
_k(\iota))=\bar{e}_1(\theta _k(\iota))=\bar{e}_1(id_E(A))=\epsilon
_1,
$$
where $\epsilon _0:E(A\otimes k[0,1])\rightarrow E(A\otimes k)$
and $\epsilon _1:E(A\otimes k[0,1])\rightarrow E(A\otimes k)$ are
homomorphisms induced by the evolutions.
\end{proof}

\begin{cor}
\label{higsr} Let $E$ be a stable and split additive functor from
an admissible subcategory of the category of real or complex
$C^{*}$-algebras to the category of abelian groups. Then it is
homotopy invariant.
\end{cor}

This result can be similarly deduced from Proposition \ref{frp} as
it is done in subsection 3.2 of \cite{Higs}.

\section {Smooth Karoubi's Conjecture for Fr\'echet $k$-algebras}

The problem of the isomorphism of algebraic and smooth
$K$-theories on the category $\mathcal{A}$ of Fr\'echet
$k$-algebras which we call Smooth Karoubi's Conjecture can be
formulate as follows.

Smooth Karoubi's Conjecture:

{\em For any Fr\'echet $k$-algebra $A$ with a bounded approximate
unit there is an isomorphism
$$
\alpha^{*}_{n}: K_{n}(A\hat\otimes \mathcal{K})
\xrightarrow{\approx} K^{sm}_{n}(A\hat\otimes \mathcal{K})
$$
for all $n\geq 0$.}

Using results of Sections 1-3 the aim of this section is to
confirm this conjecture.

To this end we will investigate the functors
$K_{n}(A\hat\otimes(-\otimes K))$, $n\in \mathbb{Z}$, on the
category of $C^{*}$-algebras, where $\otimes$ is the well-known
tensor product defined on this category. It will be shown that for
any $A$ belonging to a wide class of locally convex $k$-algebras
these functors have important homological properties such as
exactness, stability, homotopy property and Bott periodicity.

First of all the exactness property will be considered. According
to Corollary 3.12 \cite{SuWo} any locally convex $k$-algebra has
the excision property in algebraic $K$-theory if it has the
TF-property. The Cohen-Hewitt generalized theorem proved in
Section 2 allows us to establish the TF-property for a wide class
of locally convex $k$-algebras, namely for Fr\'echet $k$-algebras
with bounded approximate unit.

 Let $A$ be locally convex $k$-algebra. An element $a\in A$ is said
to be bounded with respect to a family $\{||\cdot ||_{\alpha}\}$
of seminorms if there exists a positive constant $C$ such that
$$
||a ||_{\alpha}<C.
$$

Denote by $A_b$ the space of bounded elements in $A$ with respect
to a determining family $\mathcal{F}$ of seminorms . One introduce
on $A_{b}$ a norm given by
$$
||a||= \sup_{\alpha}||a||_{\alpha},
$$
$a\in A_{b}, \alpha\in \mathcal{F}$.
\begin{thm}
{\label{higsr} Let $A$ be a Fr\'echet $k$-algebra. Then

(i) $A_{b}$ is a Banach $k$-algebra with respect to the above
defined norm and $A$ is a Fr\'echet $A_{b}$-module in the sense of
Section 2.

(ii) If $A$ is a Fr\'echet $k$-algebra with bounded approximate
unit (in the sense of Definition 2.1 as a Fr\'echet
$A_{b}$-module), then it possesses the TF-property and therefore
the excision property in algebraic $K$-theory.}
\end{thm}
\begin{proof} (i) Since $A$ is a Fr\'echet algebra, there exists a determining
countable  subset $\{||\cdot ||_n\}$ of seminorms. For any
elements $a,b \in A_b$ and a seminorm $||\cdot ||_n$ there exists
a seminorm $||\cdot ||_m$ such that
$$
||ab||_n\leq C||a||_m|\cdot||b||_m.
$$
Then one has
$$
||ab||_n\leq C||a||\cdot|||b||.
$$
Therefore $ ||ab||\leq C||a||\cdot|||b||.$ Thus $A_b$ is a normed
$k$-algebra.

Now we have to show that $A_b$ is complete. Let $(a_i)$ be a
Cauchy sequence in $A_b$. Then $(a_i)$ is a Cauchy sequence in $A$
too and it has a limit $a\in A$. If $||\cdot ||_n$ is a seminorm
and $\varepsilon >0$, there exists $i_0\in \mathbb{N}$ such that
$||a_i-a||_n<\varepsilon$ for $i\geq i_0$. Furthermore we can find
$i_1$ such that $||a_i-a_j||<\varepsilon $ for $i,j>i_1$. If we
choose $l\geq i_0,i_1$, one has
$$
||a_i-a||_n\leq ||a_i-a_l||+||a_l-a||_n\leq 2\varepsilon
\;\;\text{if}\;\;i\geq i_1
$$
and
$$
||a||_n\leq ||a-a_l||_n+||a_l-a_{i_1}||+||a_{i_1}||_n\leq
2\varepsilon +||a_{i_1}||
$$
for any $||\cdot ||_n$, implying $a\in A_b$.

It is clear that the action of $A_b$ on $A$ satisfies condition
\ref{chweqa} and that is an easy consequence of
$$
||\mu x||_n\leq C||\mu ||_m\cdot ||x||_m
$$
for some seminorm $||\cdot||_m$, where $\mu \in A_b$ and $x\in A$.

 (ii) In Theorem 2.4 replace $A$ by $A_{b}$ and $L$ by $A^{m}$. Then
for an element $x=(x_1,...,x_m)\in A^m$ there exist elements
$\sigma\in A_b$ and $y=(y_1,...,y_m)\in \overline{A_b\cdot
x}\subset \overline{A\cdot x}$ such that $x=\sigma y$. Applying
again the generalized Cohen-Hewitt factorization Theorem
\ref{chwta} to the $A_b$-module $A$, we obtain the factorization
$\sigma=\gamma\delta$ for some $\gamma,\delta\in A$ such that
$\delta \in \overline{A\gamma\delta}$. Therefore the right
annihilator $r(\delta)$ contains the right annihilator
$r(\gamma\delta)$. The inclusion $r(\delta)\subset
r(\gamma\delta)$ is always true and trivial. This proves that $A$
possesses the TF-property.
\end{proof}

Let $D$ be a real or complex $C^*$-algebra. Denote by $SB$ the
$C^*$-tensor product $k^{(0,1)}\otimes D$ and by $CB$ the
$C^*$-tensor product $k^{(0,1]}\otimes D$. If $D$ is a real
$C^*$-algebra, denote by $\mho D$ the $C^*$-tensor product
$C^{\mathbb{R}}_0(i\mathbb{R})\otimes D$, where
$C^{\mathbb{R}}_0(i\mathbb{R})$ is the real $C^*$-algebra defined
in \cite{Cunt}. Now we give some slight generalization of
Cuntz-Bott Periodicity Theorem 4.4  \cite{Cunt}, which will be
useful below.
\begin{thm}\label{analthb}Let $E:\mathcal{C}^*\rightarrow Ab$ be
a functor defined on the category $\mathcal{C}^*$ of
$C^*$-algebras and $*$-homomorphisms satisfying the following
properties:
\begin{enumerate}
    \item $E$ is homotopy invariant;
    \item $E$ is stable invariant;
    \item $E$ is half-exact that means for any proper exact sequence of $C^*$-algebras
    $$
    0\rightarrow I\rightarrow B\rightarrow C\rightarrow 0
    $$
    the sequence of abelian groups
\begin{equation}
\label{anale} E(I)\rightarrow E(B)\rightarrow E(C)
\end{equation}
is exact.
\end{enumerate}
Then there are natural isomorphisms
   \begin{itemize}
    \item $E(S^2D)\approx
    E(D)$  for any complex $C^*$-algebra $D$;
    \item $E(\mho SD)\approx
    E(D)$  for any real $C^*$-algebra $D$.
   \end{itemize}
\end{thm}

\begin{proof} Note that the exact sequences
$$
0\rightarrow \mathcal{K_{\mathbb{C}}}\rightarrow
\mathcal{T}_{\mathbb{C}}\rightarrow C_0(\mathbb{R})\rightarrow
0\;\;\;\;\text{(complex case )}
$$
and

$$
0\rightarrow \mathcal{K_{\mathbb{R}}}\rightarrow
\mathcal{T}_{\mathbb{R}}\rightarrow
\mathbb{C}_0^{\mathbb{R}}(i\mathbb{R})\rightarrow
0\;\;\;\;\text{(real case )}
$$
defined in \cite{Cunt} have bounded linear sections, since
$\mathbb{C}_0(\mathbb{R})$ and
$\mathbb{C}_0^{\mathbb{R}}(i\mathbb{R})$ are $C^*$-nuclear
$C^*$-algebras. Taking into account this observation the proof
completely coincides with the proof of the similar result in
\cite{Cunt}.
\end{proof}

We arrive to the following result.
\begin{thm}\label{analthc}
Let $A$ be a Fr\'echet $k$-algebra with a bounded approximate
unit. Then the functors
$$
K_n^{A,\mathcal{K}}=K_n(A\hat{\otimes}(-\otimes
\mathcal{K})):\mathcal{C}^*\rightarrow Ab,
$$
$n\in Z$,
\begin{enumerate}
    \item have the excision property in the following sense: if
    \begin{equation}
\label{analf}
    0\rightarrow I\rightarrow B\rightarrow C\rightarrow 0
\end{equation}
     is a proper exact sequence of $C^*$-algebras, then there is a long
    exact sequence of abelian groups
\begin{equation}
\label{anale} ...\rightarrow K_n^{A,\mathcal{K}}(C)\rightarrow
K_n^{A,\mathcal{K}}(I)\rightarrow
K_n^{A,\mathcal{K}}(B)\rightarrow
K_n^{A,\mathcal{K}}(C)\rightarrow ...
\end{equation}
    \item are stable invariant;
    \item are homotopy invariant;
    \item satisfy the following relations:
   \begin{itemize}
    \item $K_{n+1}^{A,\mathcal{K}}(D)\approx
    K_n^{A,\mathcal{K}}(SD)$ for any $C^*$-algebra $D$;
    \item $K_n^{A,\mathcal{K}}(D)\approx K_{n+1}^{A,\mathcal{K}}(SD)
    $ for any complex $C^*$-algebra $D$;
    \item $K_n^{A,\mathcal{K}}(D)\approx K_{n+1}^{A,\mathcal{K}}(\mho D)$
    for any real $C^*$-algebra $D$.
   \end{itemize}
\end{enumerate}
\end{thm}
\begin{proof}
(1). Since $\mathcal{K}$ is a $C^*$-nuclear algebra, one has the
proper exact sequence of $C^*$-algebras
$$
0\rightarrow I\otimes \mathcal{K}\rightarrow B\otimes
\mathcal{K}\rightarrow C\otimes \mathcal{K}\rightarrow 0.
$$
 According to Lemma 1.1 this implies the exactness of the following sequence of
Fr\'echet $k$-algebras
\begin{equation}
\label{analg} 0\rightarrow A\hat{\otimes}(I\otimes
\mathcal{K})\rightarrow A\hat{\otimes}(B\otimes
\mathcal{K})\rightarrow A\hat{\otimes}(C\otimes
\mathcal{K})\rightarrow 0.
\end{equation}
 Since every $C^*$-algebra has a bounded approximate unit
  and by assumption $A$ has a bounded approximate
unit too, one concludes that $A\hat{\otimes}(I\otimes
\mathcal{K})$ has also a bounded approximate unit. Thus by Theorem
4.1 the Fr\'echet $k$-algebra $A\hat{\otimes}(I\otimes
\mathcal{K})$ has the TF-property and therefore the excision
property in algebraic $K$-theory. This implies the long exact
sequence of algebraic $K$-groups associated with the short exact
sequence \ref{analg} of Fr\'echet $k$-algebras.

(2). Straightforward;

(3). (1) and (2) allow us to apply Higson's homotopy invariance
theorem (see \cite{Higs} for complex case and Corollary 3.4 for
real $C^*$-algebras) to show that the functors
$K_n^{A,\mathcal{K}}$, $n\in Z$, are homotopy invariant.

(4) Since the sequence
$$
0\rightarrow Sk\rightarrow Ck\rightarrow k\rightarrow 0
$$
 is a proper short exact sequence, the first isomorphism is
a consequence of (1) applied to the proper short exact sequence
$$
0\rightarrow A\hat{\otimes}(SB\otimes \mathcal{K})\rightarrow
A\hat{\otimes}(CB\otimes \mathcal{K})\rightarrow
A\hat{\otimes}(B\otimes \mathcal{K})\rightarrow 0.
$$
 The last two isomorphisms are immediate consequences of (1)-(3)
  and Theorem \ref{analthb}.
\end{proof}

A Fr\'echet $k$-algebra $B$ will be called {\em quasi
$\hat{\otimes}$-stable} if it has the form
$A\hat{\otimes}\mathcal{K}$ for some Fr\'echet $k$-algebra $A$
with bounded approximate unit.

 We are ready to prove the Smooth Karoubi's Conjecture.
\begin{thm}
\label{analpra} The functors $K_n(-\hat{\otimes}\mathcal{K})$ and
$K_n^{sm}(-\hat{\otimes}\mathcal{K})$, $n\geq 0$, are isomorphic
on the category of quasi $\hat{\otimes}$-stable Fr\'echet
$k$-algebras.
\end{thm}

\begin{proof}
According to Corollary \ref{smcra} it suffices to prove that the
functors  $K_n(-\hat{\otimes}\mathcal{K})$ are smooth homotopy
functors for all $n\geq 1$. Consider the commutative diagram
\begin{equation}
\label{diagr}
\begin{array}{ccc}
  K_n((A^{\infty (I^n)}\hat{\otimes }k)\hat{\otimes} \mathcal{K}) & \rightrightarrows
  & K_n((A\hat{\otimes }k)\hat{\otimes} \mathcal{K}) \\
  \downarrow \approx &  & \downarrow \approx  \\
  K_n((A\hat{\otimes }(k^{\infty (I^n)}\hat{\otimes} \mathcal{K})) & \rightrightarrows
   & K_n(A\hat{\otimes }(k\hat{\otimes} \mathcal{K}) \\
  \downarrow &  & \parallel \\
  K_n((A\hat{\otimes }(k^{I^n}\otimes \mathcal{K}))& \rightrightarrows & K_n(A\hat{\otimes }(k\otimes \mathcal{K}), \\
\end{array}
\end{equation}
where the horizontal homomorphisms are induced by evolution maps.
Since the functors $K_n(A\hat{\otimes}(-\otimes \mathcal{K}))$,
$n\in Z$, are homotopy invariant, the bottom two horizontal
homomorphisms are equal, implying the equality of the top two
horizontal homomorphisms. This shows that the functors
$K_n(-\hat{\otimes}\mathcal{K})$ are smooth homotopy invariant.
\end{proof}

\begin{thm}
\label{analthd} On the category $\overline{\mathcal{A}}$ of
Fr\'echet $k$-algebras with bounded approximate unit and
continuous $k$-homomorphisms the functors
$$
K_n^{\mathcal{K}}=K_n(-\hat{\otimes}
\mathcal{K}):\overline{\mathcal{A}}\rightarrow Ab,
$$
$n\in Z$,
\begin{enumerate}
    \item are smooth homotopy invariant;
    \item have the excision property in the following sense: if
    \begin{equation}
\label{analf}
    0\rightarrow I\rightarrow B\rightarrow C\rightarrow 0
\end{equation}
     is a proper exact sequence in $\overline{\mathcal{A}}$, then there is a long
     exact sequence of abelian groups
\begin{equation}
\label{analk} ...\rightarrow K_n^{\mathcal{K}}(C)\rightarrow
K_n^{\mathcal{K}}(I)\rightarrow K_n^{\mathcal{K}}(B)\rightarrow
K_n^{\mathcal{K}}(C)\rightarrow ...
\end{equation}
    \item $K_{n+1}^{\mathcal{K}}(A)\approx
    K_n^{\mathcal{K}}(\Omega_{sm} A)$ for any $k$-algebra $A\in Ob\overline{\mathcal{A}}$;
    \item satisfy the following relations:
   \begin{itemize}
    \item $K_n^{\mathcal{K}}(A)\approx K_{n+1}^{A,\mathcal{K}}(Sk)
    $ for any complex algebra $A\in Ob\overline{\mathcal{A}}$;
    \item $K_n^{\mathcal{K}}(A)\approx K_{n+1}^{A,\mathcal{K}}(\mho k)$ for any real algebra
    $A\in Ob\overline{\mathcal{A}}$.
   \end{itemize}
\end{enumerate}
\end{thm}

\begin{proof}
(1) is already proved for $n\geq 1$ (see proof of Theorem
\ref{analpra}) and for $n< 1$ the proof is similar.

 (2) If (4.6) is a proper short exact sequence in $\overline{\mathcal{A}}$, then so is the sequence
\begin{equation} \label{analj}
    0\rightarrow I\hat{\otimes}\mathcal{K}\rightarrow B\hat{\otimes}\mathcal{K}
    \rightarrow C\hat{\otimes}\mathcal{K}\rightarrow 0
\end{equation}
of Fr\'echet $k$-algebras implying the long exact sequence of
algebraic $K$-groups.

(3) Consider the short exact sequence \ref{smeqc}
$$
 0\rightarrow \Omega_{sm} A\rightarrow
\mathcal{I}(A)\xrightarrow{\tau _A} A\rightarrow 0,
$$
where the epimorphism $\tau _A$ has a natural bounded section
given by $a\mapsto a\cdot t$, $t\in [0,1]$. Then the exact
sequence
\begin{equation}
 0\rightarrow \Omega_{sm} A\hat{\otimes}\mathcal{K}\rightarrow
\mathcal{I}(A)\hat{\otimes}\mathcal{K}\xrightarrow{\tau _A}
A\hat{\otimes}\mathcal{K}\rightarrow 0
\end{equation}
is a proper short exact sequence. The Fr\'echet $k$-algebra
$\Omega_{sm} A\hat{\otimes}\mathcal{K}$ has a bounded approximate
unit and therefore (4.9) induces the long exact sequence of
algebraic $K$ -groups. Since
$K_n(\mathcal{I}(A)\hat{\otimes}\mathcal{K})=0$, one gets the
required isomorphism $K_{n+1}^{\mathcal{K}}(A)\approx
    K_n^{\mathcal{K}}(\Omega_{sm} A)$.

(4) is an immediate consequence of (4) of Theorem \ref{analthc}.
\end{proof}

\begin{rem}
There is another way to prove Theorem \ref{analpra}. Since
$K_n^{sm}(-\hat\otimes \mathcal{K})$, $n\geq 1$, are smooth
homotopy functors, one has isomorphisms
$$
K_n^{sm}(A\hat{\otimes }\mathcal{K})\approx K_1^{sm}(\Omega_{sm}
^{n-1}A\hat{\otimes }\mathcal{K}),\;\; n\geq 1,
$$
and $K_1^{sm}(\Omega_{sm} ^{n-1}A\hat{\otimes }\mathcal{K})\approx
K_1(\Omega _{sm}^{n-1}A\hat{\otimes }\mathcal{K})$ by Theorem
\ref{smtha}. On the other hand, $K_1(\Omega_{sm}
^{n-1}A\hat{\otimes }\mathcal{K})\approx K_n(A\hat{\otimes
}\mathcal{K})$ by Theorem \ref{analthd}(3). Therefore
$$
K_n^{sm}(A\hat{\otimes }\mathcal{K})\approx K_1(\Omega_{sm}
^{n-1}A\hat{\otimes}\mathcal{K})\approx
K_n(A\hat{\otimes}\mathcal{K}),
$$
$n\geq 1$.
\end{rem}

\begin{rem}
A locally convex $k$-algebra $A$ has the smooth homotopy (the
homotopy) property in algebraic $K$-theory if there is an
isomorphism $K_{*}(A)\xrightarrow {\approx} K_{*}(A^{\infty
(I^{n})})$ ($K_{*}(A)\xrightarrow {\approx}K_{*}(A^{(I^{n})})$)
for all $n> 0$. It seems that "the homotopy property in algebraic
K-theory" is more appropriate than the corresponding
"$K_{*}$-stability" used in \cite{Inas1}, Definition 12. By
Theorem 4.5(1) any quasi $\hat\otimes$-stable Fr\'echet
$k$-algebra has the smooth homotopy property in algebraic K-theory
and it is clear that any stable $C^{*}$-algebra has the homotopy
property in algebraic K-theory.
\end{rem}
\

\begin{center}
{\bf Acknowledgements}
\end{center}

\;

The authors were partially supported by INTAS grant 00 566, INTAS
grant 03-51-3251 and GRDF grant GEM1-3330-TB-03.

\end{document}